\documentclass[12pt]{article}
\newif\ifhide
\hidefalse

\usepackage{amssymb}
\date{}

\title{%{\normalsize\tt\hfill\jobname.tex}\\
On large deviations in the averaging principle for SDE's with a
``full dependence'', correction
\author{A. Yu. Veretennikov
\\
{(School of Mathematics, University of Leeds, UK} \\
{\& Institute for Information Transmission Problems, Russia)} } }

\begin{document}
\maketitle
\newcommand{\eps}{\varepsilon}
\newcommand{\argmax}{\mathop{\rm argmax}\nolimits}
\newtheorem{Theorem}{Theorem}
\newtheorem{Lemma}{Lemma}

{\footnotesize
We establish the large deviation principle for stochastic
differential equations with averaging in the case when all
coefficients of the fast component depend on the slow one,
including diffusion.}
\footnote{{\it AMS} 1991 {\it subject classifications.}
60F10, 60J60.}
\footnote{{\it Key words and phrases.} Large deviations, averaging,
stochastic differential equation.}

%\footnotemark{Key words and phrases: large deviations, averaging,
%stochastic differential equation
%AMS 1991 subject classifications. 60F10, 60J60.}

\section{Introduction}
This is a corrected version of the paper \cite{Veretennikov
(1999)}. We
consider the SDE system
\begin{eqnarray}\label{eqAA}
& dX_t = f(X_t,Y_t)dt, \quad X_0=x_0, \nonumber
 \\
& dY_t = \eps^{-2}B(X_t,Y_t)dt + \eps^{-1}C(X_t,Y_t)dW_t, \quad
Y_0=y_0.
\end{eqnarray}
Here $X_{t}\in E^{d}, \, Y_{t}\in M$, $M$ is a compact manifold
of dimension $\ell$ (e.g. torus $T^\ell$), $f$ is a function
with values in $d$--dimensional Euclidean space $E^{d}$, $B$ is
a function with values in $T M$, $C$ is a function with values
in $(T M)^\ell$ (i.e., in local coordinates an $\ell\times \ell$
matrix), $(W_{t})$ is an $\ell$--dimensional Wiener process
with respect to some increasing and right continuous
filtration $({\cal F}_t)$ on some
probability space $(\Omega, F, P)$, $\eps>0$ is a small
parameter, i.e., $\eps \to 0$.  Concerning SDE's on
manifolds we refer to \cite{Ikeda and Watanabe (1989)}.

The large deviation principle (LDP) for such systems with a ``full
dependence'', that is, $C(X_t,Y_t)$, was not treated before
\cite{Veretennikov (1999)}. Only
the case $C(Y_t)$ was considered in the papers by
\cite{Freidlin (1976),
Freidlin (1978), Freidlin and Wentzell (1984)} for a compact
state space and by \cite{Veretennikov (1994)} for a
non-compact one. Also the papers \cite{Liptser (1996)}, 
\cite{Veretennikov (1998)} and \cite{Liptser (2002)} on
similar or close topics for more general systems with small
additive diffusions should be mentioned, which, however,
all concern only the case $C(Y_t)$. Concerning most recent
developments the reader is referred to \cite{Kifer} and the
references therein.

The LDP for systems like (\ref{eqAA}) is important in averaging
and homogenization, in the KPP equation theory, for stochastic
approximation algorithms with averaging and so forth. The problem
of an LDP for the case $C(X_t,Y_t)$ has arisen since
\cite{Freidlin (1976), Freidlin (1978)}. Intuitively, the
scheme used for $C(Y_t)$ should work; at least, almost all
main steps go well. Indeed, there was only one lacuna; the
use of Girsanov's transformation did not allow freezing of
$X_t$ if $C$ depended on the slow motion, while it worked
well and very naturally for the drift $B(X_t,Y_t)$. Yet the
problem remained unresolved for years and
the answer was unclear. Notice that this difficulty does
not appear in analogous discrete-time systems (see
\cite[Chapter 11]{Gulinsky and Veretennikov (1993)}).

It turned out that the use of Girsanov's transformation in some
sense prevented from resolving the problem. Our approach in
this paper is based on a new technical result, Lemma
\ref{L-main} below. The main new idea is to use two
different scales of partitions of the interval $[0,T]$, a
``first-order
partition'' by points $\Delta,\, 2\Delta,\, \ldots$, which do not
depend on the small parameter $\eps$ and ``second-order
partitions'' which depend on $\eps$ in a special way, by points
$\eps^2 t(\eps),\, 2\eps^2 t(\eps), \ldots\,$.  Then
the exponential estimates needed for the proof of the
result can be established in two steps. First, the estimates
for a ``small'' partition interval are derived
using the uniform bound of Lemma \ref{L-lim} (see
below) and the estimates for stochastic integrals. It is
important that in the ``second'' scale the fast motion is still close
enough to its frozen version [the bound (\ref{eqYY}) below].
Second, the bounds for ``small'' partitions and induction give
one the estimate for a ``large'' partition interval.

{\it The original proof in \cite{Veretennikov (1999)}
contained some gap 
relates to a boundedness of some auxiliary constant $b$ in
the proof: in the original version this constant may
depend implicitly on the partition size $\Delta$, while the choice of $\Delta$ 
could depend on $b$, hence generating a vicious circle.
The main aim of this version of the paper is to present the
``patch''. A provisional version of this correction may be
found in \cite{Veretennikov 2005}. The present version is
simplified further. The correction uses improved
approximations that keep this constant $b$ bounded in the
lower and upper bounds, and it uses also a truncated
Legendre transformation in the upper bound.
The author is deeply indebted to Professor Yuri Kifer for
discovering this vicious circle in the original version of the
paper. The main technical tool remains the Lemma \ref{L-main}. All
standing assumptions are the same as in the original version.}

The main result is stated in Section 2. In Section 3 we
present
auxiliary lemmas, among them the main technical Lemma
\ref{L-main}
with its proof and a version of an important lemma from
\cite{Freidlin and Wentzell (1984)} (see Lemma \ref{L-FW})
which requires certain
comments. Those comments along with other related remarks
are
given in the Appendix, the latter has been also slightly extended.
The proof of the main theorem is presented in Section 4.

\section{Main result}
We make the following assumptions.

\begin{description}
\item[$(A_f)$]
The function $f$ is bounded and satisfies the Lipschitz condition.

\item[$(A_C)$]
The function $CC^*$ is bounded, uniformly nondegenerate, $C$
satisfies the Lipschitz condition.

\item[$(A_B)$]
The function $B$ is bounded and satisfies the Lipschitz
condition.
\end{description}

Some conditions may be relaxed; for example, $B$ may be
assumed locally bounded, $C$ locally (with respect
to  $x$)
nondegenerate and so on.

The family of processes $X^\eps$ satisfies a large deviation
principle in the space $C([0,T];R^d)$ with a normalizing
coefficient $\eps^{-2}$ and a rate function $S(\varphi)$ if
the following three conditions are satisfied:
\begin{equation}\label{eq4}
\limsup_{\eps\to 0} \eps^2 \log P_x(X^\eps \in F) \leq - \inf_F
S(\varphi), \quad \forall F \mbox{ closed },
\end{equation}
\begin{equation}\label{eq5}
\liminf_{\eps\to 0} \eps^2 \log P_x(X^\eps\in G) \geq - \inf_G
S(\varphi), \quad \forall G \mbox{ open },
\end{equation}
and $S$ is a ``good'' rate function; that is, for any $s\geq 0$,
the set
\begin{displaymath}
\Phi(s) : = (\varphi\in C([0,T];R^d):
\,\, S(\varphi)\leq s, \,\, \varphi(0)=x)
\end{displaymath}
is compact in $C([0,T];R^d)$. We will establish the
following equivalent set of assertions due to Freidlin and
Wentzell, where $\rho(\phi, \psi) = \sup_{0\le s\le
T}|\phi_s - \psi_s|$, 
\begin{equation}\label{eq4a}
%\limsup_{\delta\to 0}
\limsup_{\delta\to 0} \limsup_{\eps\to 0} \eps^2 \log
P_x(\rho(X^\eps, \Phi(s))\ge \delta) \leq - s, \quad
\forall
s>0,
\end{equation}
where $\Phi(s) := \{\varphi\in C[0,T;R^d], \, S(\varphi)\le s\}$,
and
\begin{equation}\label{eq5a}
\liminf_{\delta\to 0}\liminf_{\eps\to 0} \eps^2 \log
P_x(\rho(X^\eps,\varphi)<\delta) \geq - S(\varphi),
\quad
\forall \varphi,
\end{equation}
where $S$ is a ``good'' rate function (see above). In what
follows, $\dot\varphi_t$ is a derivative function for
$\varphi_t$ and if it does not exist almost everywhere or
if the integral $\int_0^T L(\varphi_t,\dot\varphi_t)\,dt$
diverges, then by definition $\int_0^T
L(\varphi_t,\dot\varphi_t)\,dt := +\infty$.

\begin{Theorem}\label{Theorem1}
Let $(A_f)$, $(A_B)$, $A_C)$ be satisfied. Then the family
$(X^\eps_t = X_t,\; 0\leq t\leq T)$ satisfies the LDP as $\eps
\to 0$ in the space $C([0,T];R^d)$ with a rate function
$$
  S(\varphi) = \int_0^T L(\varphi_t,\dot\varphi_t)\,dt,
$$
where
$$
  L(x,\alpha) = \sup_\beta (\alpha\beta - H(x,\beta)),
$$
$$
  H(x,\beta) = \lim_{t\to\infty} t^{-1} \log
E\exp\left(\beta \int_0^t f(x,y^x_s)ds\right).
$$
The limit $H$ exists and is finite for any $\beta$, the
functions $H$ and $L$ are convex in their last arguments $\beta$
and $\alpha$ correspondingly, $L\geq 0$ and $H$ is continuously
differentiable in $\beta$.
\end{Theorem}

The differentiability of $H$ at any $\beta$ will be provided
by the compactness of the state space of the fast
component. The constants $C$ in the calculus may change
from line to line, unlike $K, C_f, L_f$ and some other.

\section{Auxiliary lemmas}
Let $\tilde W_t = \eps^{-1} W_{t\eps^{2}}$, $y_t =
Y_{t\eps^{2}}$, $x_t = X_{t\eps^{2}}$, and let $y^x_t$
solve an SDE,
\begin{equation}\label{eqBB}
dy^x_t = B(x,y^x_t)dt + C(x,y^x_t)d\tilde W_t, \quad
y^x_0=y_0.
\end{equation}
%and $Y^x_t := y^x_{t/\eps^2}$.
Below $\tilde {\cal F}_t: = {\cal F}_{t\eps^2}$,
$\beta\in E^d$, $\beta f$ means a scalar product and the
index $y$ in $E_y$ stands for the initial value of $y_t^x$
at $t=0$. 
Let us consider the semigroup of operators $T^{\beta }_{t},
t\geq
0$,  on $C(M)$ defined by the formula
$$
T^{x',x,\beta}_{t}g(y)= T^{\beta }_{t}g(y) =
E_{y}g(y^{x}_{t})\exp \left(\int^{t}_{0}
\beta f(x',y^{x}_{s})ds\right),
$$
where $\beta\in E^d$, $\beta f$ is a scalar product and the
index $y$ in $E_y$ means the initial value of $y_t^x$ at
$t=0$. In the case if some inequality is uniform over $y\in
M$, this index may be dropped in the calculus.

\begin{Lemma} \label{L-comp}%Le1
Let assumptions $(A_f)$, $(A_B)$, $(A_C)$ be satisfied.
Then for any $\beta$ the operator $T^{\beta }_{1}$ is compact in
the space $C(M)$.
%$C([0,T];R^d)$.
\end{Lemma}

\begin{Lemma} \label{L-kls}%Le2
Let assumptions $(A_f)$, $(A_B)$, $(A_C)$ be satisfied. Then
the spectral radius $r(T^{\beta }_{1})$ is a simple eigenvalue
of $T^{\beta }_{1}$ separated from the rest of the spectrum
and its eigen--function $e_{\beta }$ belongs to the cone
$C^+(M)$. Moreover, function $r(T^{\beta }_{1})$ is smooth (of
$C^\infty$) in $\beta$ and for any $b>0$ the function
$e_{\beta }$ is bounded and separated away from zero
uniformly in $|\beta|<b$ and all $x', x$.
\end{Lemma}

\begin{Lemma} \label{L-lim}%Le3
Let $\beta \in E^{d}$, and let assumptions $(A_f)$, $(A_B)$,
$(A_C)$ be satisfied. Then there exists a limit uniformly
in $x,x'$,
$$
\tilde H(x',x,\beta ) = \lim_{t\to \infty} t^{-1}\log
E_{y}\exp\left(\beta\int
^{t}_{0} f(x',y^{x}_{s})ds\right);
$$
moreover, $\tilde H(x',x,\beta ) = \log
r(T^{x',x,\beta}_1)$.
The function $\tilde H(x',x,\beta )$ is of $C^\infty$ in
$\beta $ and convex in $\beta $. For any $b>0$ there
exists $C(b)$ such that, for any $y$, $|\beta|<b$ and for
all values of $t>0$ uniformly in $x,x'$,
\begin{equation}\label{eqUU}
|t^{-1}\log E_{y}\exp \left(\beta \int^t_0
f(x',y^{x}_{s})ds\right) - \tilde H(x',x,\beta )|\le
C(b)t^{-1}.
\end{equation}
Notice that $|\tilde H(x',x,\beta)|\leq \|f\|_C |\beta|$.
\end{Lemma}
In what follows, $\nabla_\beta \tilde H$ stands for the
gradient of $\tilde H$ with respect to $\beta$.
\begin{Lemma} \label{L-dif}%Le4
Let assumptions $(A_f)$, $(A_B)$, $(A_C)$ be satisfied. Then
for any $b>0$ the functions $\tilde H$ and $\nabla_\beta \tilde H$ 
are uniformly continuous in $(x',x,\beta )$, for $|\beta |<b$.
\end{Lemma}

Lemmas \ref{L-comp} -- \ref{L-dif} are standard (cf.
\cite{Veretennikov (1994)} or \cite{Veretennikov (1992)}).
They are based on Frobenius-type theorems for positive
compact operators (see \cite{Krasnoselskii (1989)})
%\cite{Krasnosel'skii, Lifshitz and Sobolev (1989)}) 
and the theory of perturbations of linear operators (see
\cite[Chapter 2]{Kato (1976)}). 
%Denote $\tilde F_t = F_{t\eps^{2}}$.

\begin{Lemma}\label{L-main}
Let the assumptions $(A_f)$, $(A_B)$, $(A_C)$ hold true,
$b>0$, $t(\eps) \to \infty$ and $t(\eps) = o(\log
\eps^{-1})$ as $\eps \to 0$. Then
for any $\nu>0$ there exist $\delta(\nu)>0$, $\eps(\nu)>0$
such
that for $\eps \leq \eps(\nu)$ uniformly with respect
to 
$t_0,\,x',\,x,\,x_{0},\,y_{0}$, $x_{t_0}$, $|\beta|\leq
b$,
 the inequality holds on the set
$\{|x_{t_{0}} - x| < \delta(\nu)\}$,

\begin{equation}\label{eqCC}
  \left|\log E_{}(\exp(\beta
\int\limits_{t_{0}}^{t_{0}+t(\eps)}
f(x',y_s)ds)|\tilde {\cal F}_{t_{0}})
  - t(\eps) \tilde H(x',x,\beta)\right| \leq \nu t(\eps).
\end{equation}
Moreover, if $\Delta \leq \Delta(\nu) =
(1+\|f\|_C)^{-1}\delta(\nu)/2$ and $\eps$ is small enough, then
uniformly with respect
to  $t_0,\,x',\,x,\,x_{0},\,y_{0}$,
$\delta\le\delta(\nu)$, $|x_{t_0}-x|<\delta$, and
$|\beta|\leq b$,
\begin{eqnarray}\label{eqDD}
& \exp(\eps^{-2}\Delta \tilde H(x',x,\beta) -
\nu\Delta\eps^{-2})
  \nonumber \\
&\leq E_{}\left( \exp(\beta\eps^{-2}
\int\limits_{t_{0}}^{t_{0}+\Delta} f(x',Y_s)ds ) |
{\cal F}_{t_{0}}\right)
  \nonumber \\
&\leq \exp(\eps^{-2} \Delta \tilde H(x',x,\beta) + \nu\Delta\eps^{-2}).
\end{eqnarray}
\end{Lemma}
{\bf Remark.} Let us emphasize that any couple
$(\Delta,\delta)$ satisfying only $\Delta \le \Delta(\nu)$ and
$\delta \le \delta(\nu)$ would do.

\noindent {\bf Proof.} {\it Step} 1. It suffices to
prove (\ref{eqCC}) and (\ref{eqDD}) for $t_0=0$. Moreover,
since $H$ is continuous, it suffices to check both
inequalities for $x=x_0$.
Indeed, the bound
$$
\left|\log E_{} \exp\left(\beta \int\limits_{0}^{t(\eps)}
f(x',y_s)ds\right) - t(\eps) \tilde H(x',x_0,\beta)\right| \leq \nu
t(\eps)
$$
implies
\begin{eqnarray*}
&\left|\log E_{} \exp\left(\beta \int\limits_{0}^{t(\eps)}
f(x',y_s)ds\right) - t(\eps) \tilde H(x',x,\beta)\right|
 \\
&\leq t(\eps)(\nu + |\tilde H(x',x,\beta)-\tilde H(x',x_0,\beta)|),
\end{eqnarray*}
and we use the uniform continuity of the function $H$ on compact
sets (remind that $|\beta|\le b$).  The same arguments are
applicable to the second inequality of the assertion of the lemma.
So, in the sequel we consider the case $x_0=x$.

Let us show first that
\begin{equation}\label{eqEE}
\sup_{x',x_{0}}\left| t(\eps)^{-1} \log
E_{}\exp\left(\beta
\int_0^{t(\eps)} f(x',y_s)ds\right) - \tilde H(x',x,\beta)\right|
\leq \nu
\end{equation}
if $\eps$ is small enough. Due to Lemma \ref{L-lim}, it would be
correct if $y_s$ were replaced by $y^x_s$ and $t(\eps)\geq
\nu^{-1}C(b)$. We will also use the bounds
\begin{equation}\label{eqXX}
\sup_{0\leq s\leq t} |x_s-x_0| \leq \eps^2 t \|f\|_C, \quad
\& \quad
\exp(Ct(\eps))t(\eps)^2 \eps^2 \to 0 \;\; (\forall C), \;\;
\eps\to 0.
\end{equation}
Let $|f(x',y)-f(x',y')|\leq L_f |y-y'|$ for all $y,y',x'$,
$L_f>0$, $C_f=\|f\|_C$. We estimate for $t(\eps) > \nu^{-1}
C(b)/4$,
\begin{eqnarray}\label{eqFF}
E_{}\exp\left(\beta\int\limits_0^{t(\eps)}
f(x',y_s)ds\right) \hspace{4cm}
 \nonumber\\
\times
\left\{I\left(\sup_{0\leq t\leq t(\eps)} |y_t-y^x_t|
\leq \nu/(4L_fb)\right)
+ I\left(\sup_{0\leq t\leq t(\eps)} |y_t-y^x_t| > \nu/(4L_fb)
\right)\right\}
  \nonumber \\
\leq E_{}\exp\left(\beta\int\limits_0^{t(\eps)}
f(x',y^x_s)ds +
t(\eps)\nu/4\right)
I\left(\sup_{0\leq t\leq t(\eps)} |y_t-y^x_t|\leq \nu/(4L_fb)
\right)
  \nonumber \\
+ \exp(C_f b t(\eps)\nu) E_{} I\left(\sup_{0\leq t\leq
t(\eps)}
|y_t-y^x_t| > \nu/(4L_fb)\right) \hspace{2cm}
  \nonumber \\
\leq E_{}\exp\left(\beta\int\limits_0^{t(\eps)}
f(x',y^x_s)ds\right)
\exp(t(\eps)\nu/4) 
  \hspace{3cm}
  \nonumber \\
+ 16 L_f^2 b^2\,\exp\left(C_f b t(\eps)\nu\right)\nu^{-2}
E_{}\sup_{t\leq t(\eps)}|y_t-y^x_t|^2. 
   \hspace{2,5cm}
\end{eqnarray}
By virtue of the Lemma \ref{L-lim} we have
\begin{equation}\label{eqaux1}
E_{}\exp\left(\beta\int\limits_0^{t(\eps)}
f(x',y^x_s)ds\right)
\leq \exp(t(\eps)(\tilde H(x',x,\beta)+\nu/4)),
\end{equation}
if $\eps$ is small enough.

~

Let us estimate the second term in (\ref{eqFF}). By virtue of the inequalities
for the It\^o and Lebesgue integrals, we have
\begin{eqnarray*}
E_{}\sup_{t'\leq t}|y_{t'}-y^x_{t'}|^2 %\hspace{3cm}
 %\\
\leq C E_{} \int\limits_0^t |C(x_s,y_s) - C(x,y^x_s))|^2
ds
 \\
+ Ct E_{}\int\limits_0^t |B(x_s,y_s) - B(x,y^x_s))|^2 ds
 \\
 \leq C\int\limits_0^t E|x_s-x|^2ds
+ C\int\limits_0^t E_{}\sup_{u\leq s} |y_s-y^x_s|^2ds
 \\
\leq Ct^2\eps^2 + C\int\limits_0^t E_{}\sup_{u\leq s}
|y_u-y^x_u|^2
ds.
\end{eqnarray*}
By virtue of Gronwall's lemma, one gets
$$
E_{}\sup_{t'\leq t} |y_{t'}-y^x_{t'}|^2 \leq
Ct^2\eps^2\exp(Ct).
$$
In particular,
\begin{equation}\label{eqYY}
E_{}\sup_{t'\leq t(\eps)} |y_{t'}-y^x_{t'}|^2 \leq
Ct(\eps)^2\eps^2\exp(Ct(\eps)).
\end{equation}
So the second term in (\ref{eqFF}) does not exceed the value
$\exp(C_f b \,t(\eps)\nu)\nu^{-2}Ct(\eps)^2\eps^2$ which is
$o(\exp(K t(\eps)))$ for any $K<0$. Indeed, for any such $K$ we have, 
$\exp(t(\eps)(C_fb\nu-K))\nu^{-2}Ct(\eps)^2\eps^{2} \to 0$, 
because $\exp(t(\eps)C)$ for any $C>0$ increases slower than $\eps^{-2}$ 
due to the assumption $t(\eps) = o(\log \eps^{-1}), \; \eps \to 0$. 
Hence, we get with any $K<0$ for $\eps>0$ small enough,
\begin{eqnarray*}
E_{}\exp\left(\beta\int\limits_0^{t(\eps)}
f(x',y_s)ds\right) \hspace{3cm}
 \\\\
\leq E_{}\exp\left(\beta\int\limits_0^{t(\eps)}
f(x',y^x_s)ds\right)
\exp(t(\eps)\nu/4)  + C\exp\left(K t(\eps)\right)
 \\\\
\le \exp(t(\eps)(\tilde H(x',x,\beta)+\nu/2))  
+ C\exp\left(K t(\eps)\right), 
\end{eqnarray*}
by virtue of (\ref{eqaux1}). The upper bound in (\ref{eqEE}) follows.

~

The lower bound in (\ref{eqEE}) may be etablished similarly.  
For the convenience of the reader we show the calculus. 
We estimate for $t(\eps) > \nu^{-1} C(b)/4$,
\begin{eqnarray}\label{eqFFlow}
E_{}\exp\left(\beta\int\limits_0^{t(\eps)}
f(x',y^{x}_s)ds\right) \hspace{4cm}
 \nonumber\\
\times
\left\{I\left(\sup_{0\leq t\leq t(\eps)} |y_t-y^x_t|
\leq \nu/(4L_fb)\right)
+ I\left(\sup_{0\leq t\leq t(\eps)} |y_t-y^x_t| > \nu/(4L_fb)
\right)\right\}
  \nonumber \\
\leq E_{}\exp\left(\beta\int\limits_0^{t(\eps)}
f(x',y_s)ds + t(\eps)\nu/4\right)
I\left(\sup_{0\leq t\leq t(\eps)} |y_t-y^x_t|\leq \nu/(4L_fb)
\right)
  \nonumber \\
+ \exp(C_f b t(\eps)\nu) E_{} I\left(\sup_{0\leq t\leq
t(\eps)}
|y_t-y^x_t| > \nu/(4L_fb)\right) \hspace{2cm}
  \nonumber \\
\leq E_{}\exp\left(\beta\int\limits_0^{t(\eps)}
f(x',y_s)ds\right)
\exp(t(\eps)\nu/4) 
  \hspace{3cm}
  \nonumber \\
+ 16 L_f^2 b^2\,\exp\left(C_f b t(\eps)\nu\right)\nu^{-2}
E_{}\sup_{t\leq t(\eps)}|y_t-y^x_t|^2. 
   \hspace{2,5cm}
\end{eqnarray}
Since the second term in (\ref{eqFF}) is
$o(\exp(K t(\eps)))$ with any $K<0$, this implies the bound
\begin{eqnarray*}
E_{}\exp\left(\beta\int\limits_0^{t(\eps)}
f(x',y^{x}_s)ds\right) \hspace{3cm}
 \\\\
\leq E_{}\exp\left(\beta\int\limits_0^{t(\eps)}
f(x',y_s)ds\right)
\exp(t(\eps)\nu/4) + C\exp\left(K t(\eps)\right), 
\end{eqnarray*}
or, equivalently, 
\begin{eqnarray*}
E_{}\exp\left(\beta\int\limits_0^{t(\eps)}
f(x',y_s)ds\right)\hspace{3cm}
 \\\\
\ge 
E_{}\exp\left(\beta\int\limits_0^{t(\eps)}
f(x',y^{x}_s)ds\right) 
\exp(-t(\eps)\nu/4) - C\exp\left(K t(\eps)\right). 
\end{eqnarray*}
Now due to (\ref{eqaux1}), we get, with any $K<0$ and 
$\eps>0$ small enough,
\begin{eqnarray*}
E_{}\exp\left(\beta\int\limits_0^{t(\eps)}
f(x',y_s)ds\right)
% \\\\
\ge \exp(t(\eps)(\tilde H(x',x,\beta)-\nu/2)) 
- C\exp\left(K t(\eps)\right),  
\end{eqnarray*}
which implies the lower bound in (\ref{eqEE}). 

~

Notice that both bounds in (\ref{eqEE}) are uniform with respect to 
$|\beta|\leq b$ and $x',\, x,\, y_0$. Since the function $H$ is
continuous, we get on the set $\{|x_{t_{0}} - x| <
\delta(\nu)\}$,
\begin{eqnarray}\label{eqG}
\sup\limits_{x',x,y_{0},t_{0},|\beta|\leq b} \left| \log
E_{}\left(\exp\left(\beta
\int\limits_{t_{0}}^{t_{0}+t(\eps)}
f(x',y^x_s)ds\right) \mid \tilde  {\cal F}_{t_{0}}\right)
%\right. \nonumber \\ \nonumber \\ \left.
- t(\eps) \tilde H(x',x,\beta)\right| \leq \nu t(\eps)
\end{eqnarray} if
$\delta(\nu)$ is small enough. 

~

\noindent
{\it Step} 2.
Let $\Delta \leq (1+\|f\|_C)^{-1} \delta (\nu)/2 = \Delta(\nu)$
and $N=\Delta
\eps^{-2} t(\eps)^{-1}$. Then $\sup_{0\leq s\leq Nt(\eps)} |x_s
- x_0| \leq \delta(\nu)/2$. Let $|x-x_0|<\delta(\nu)/2$. So,
$\sup_{0\leq s\leq Nt(\eps)} |x_s - x| < \delta(\nu)$. In
particular, $|x_{k t(\eps)} - x| < \delta(\nu)$ for any $1 \leq
k \leq N$. By induction, we get from (\ref{eqG}) for such $k$,

\begin{eqnarray*}
\exp(k t(\eps) \tilde H(x',x,\beta) - \nu k t(\eps))
 \\
\leq E_{} \exp\left(\beta \int\limits_0^{k t(\eps)}
f(x',y_s)ds
\right)
 \\
\leq \exp(k t(\eps) \tilde H(x',x,\beta) + \nu k t(\eps) ),
\end{eqnarray*}
or, after the time change,
\begin{eqnarray*}
\exp(k t(\eps) \tilde H(x',x,\beta) - \nu k t(\eps))
 \\
\leq E_{} \exp\left(\beta\eps^{-2}
\int\limits_0^{k t(\eps)\eps^{-2}}
f(x',Y_s)ds\right)
 \\
\leq \exp(k t(\eps) \tilde H(x',x,\beta)
+\nu k t(\eps)).
\end{eqnarray*}
Since $H$ is continuous then we obtain for $k=N$,
\begin{eqnarray}\label{eq8}
\exp(\eps^{-2}\Delta \tilde H(x',x,\beta) -
\nu\Delta\eps^{-2})
  \nonumber \\ \nonumber \\
\leq E_{} \exp\left(\beta\eps^{-2}
\int\limits_0^\Delta f(x',Y_s)ds\right)
 \nonumber \\ \nonumber \\
\leq \exp(\eps^{-2}\Delta \tilde H(x',x_0,\beta)
+ \nu\Delta\eps^{-2}).
\end{eqnarray}
The Lemma \ref{L-main} is proved. \hfill QED

%In the sequel we denote $\Delta = \Delta(\nu)$.
%An important point is that $\delta(\nu)/\Delta(\nu) = const > 0$.
The next Lemma is an improved version of the Lemma 7.5.2
from \cite{Freidlin and Wentzell (1984)}. Although we will
not use it explctly, its technique is essential.

\begin{Lemma}\label{L-FW}(\cite{Freidlin (1978), Freidlin
and Wentzell (1984)}). Let $S(\varphi)<\infty$. If $\psi^n$
is a
sequence of step functions tending uniformly to $\varphi$ in
$C[0,T];R^d)$ as $n\to\infty$, then there exists a sequence of
piecewise linear functions $\chi^n$ (with the same partitions)
which also tend uniformly to $\varphi$ and such that
$$
  \limsup_{n\to\infty} \int_0^T L(\psi^n_s,
\dot\chi^n_s) ds\leq S(\varphi).
$$
Moreover, one may assume without loss of generality that for any
$s$ there exists a value
$$
\beta_s = \argmax\limits_{\beta} (\beta \dot\chi^n_{s+} -
\tilde H(\psi^n_s,\psi^n_s,\beta))
$$
and
$$
L(\psi^n_s,\alpha) > L(\psi^n_s,\dot\chi^n_{s+}) +
(\alpha-\dot\chi^n_{s+})\beta_s \quad \forall
\alpha\not=\dot\chi^n_s.
$$
If $\hat\psi$ is close enough to $\psi^n_s$ then there exists
a value
$$
\hat\beta_s = \argmax\limits_{\beta} (\beta \dot\chi^n_{s+} -
\tilde H(\psi^n_s,\hat\psi,\beta)),
$$
$$
L(\psi^n_s,\hat\psi,\alpha)
> L(\psi^n_s,\hat\psi,\dot\chi^n_{s+}) +
(\alpha-\dot\chi^n_{s+})\hat\beta_s \quad \forall
\alpha\not=\dot\chi^n_s
$$
and
$$
L(\psi^n_s,\hat\psi,\dot\chi^n_{s+}) \to
L(\psi^n_s,\psi^n_s,\dot\chi^n_{s+}), \quad \hat\psi\to\psi^n_s.
$$
\end{Lemma}

We added to the original assertion the property
that $\chi^n_t$ may be chosen piecewise
linear. Indeed, such functions are used in the proof; see
\cite[Section 7.5]{Freidlin and Wentzell (1984)}.
The existence of $\beta_s$ asserted in the lemma also follows
from the proof; see \cite{Freidlin (1978)} or \cite{Freidlin
and Wentzell (1984)}.
Assertions about $\hat\psi$ and $\hat\beta_s$ also added to
the original assertion can be deduced from the proof using
similar arguments.

In fact, there is a little gap in the original proof, namely, an
additional assumption was used which was not formulated
explicitly. This is why we present a precise
statement and
give necessary comments in the Appendix.

%~

\section{Proof of theorem 1} %: non-degenerate case}
%\begin{enumerate}
%\item
{\bf 1.} 
\underline{\bf First part of the proof:}  the lower bound. Let
$S(\varphi)<\infty$, and $\nu > 0$. To establish the lower bound,
we will show the inequality: given any $\nu>0$, and {\it any}
$\delta>0$, we have for $\eps
>0$ small enough,
$$% \begin{equation}\label{eq5a}\liminf_{\delta\to 0}\liminf_{\eps\to 0}
\eps^2 \log P_x(\rho(X^\eps,\varphi)<\delta) \geq -
S(\varphi) - \nu.
$$
%\item %Step 1
Denote $H(x,\beta) = \tilde H(x,x,\beta)$. The existence of
the limit
$ \tilde H(x,x',\cdot)$ for any $x, x'$, and its
differentiability and
continuity are asserted in Lemmas 3 and 4. Throughout the
proof, we may and will assume that for any $s$,
$L(\varphi_s,\dot\varphi_s)<\infty$. Indeed, this may be
violated
only on a set of $s$ of Lebesgue measure zero. Notice that due to
the boundedness of the function $f$, this inequality implies
$\sup_s|\dot\varphi_s|\le \|f\|_C$, since for any $|\alpha|>
\|f\|_C$, we have $L(x,\alpha)=+\infty$.
Unlike in the previous section, in the sequel
both $X_0=x_0$ and $Y_0=y_0$ are fixed, hence, the symbols
$P$ and $E$ will be used without indices.

~

%\item 
\noindent
{\bf 2.}
We are going to reduce the problem of estimation from below the
probability
$$
P(\rho(X,\varphi)<\delta)
$$
to that for the probability
$$
P(\rho(X^\varphi,\varphi)<\delta'), \quad  \mbox{where} \quad
X^\psi_t := x_0 + \int_0^t f(\psi_s,Y_s)ds, \; \forall \psi,
$$
and further to
$$
P(\rho(X^{\psi},\chi)<\delta'),
$$
where both $\psi,\chi$ approximate $\varphi$. The rough idea is
eventually to choose a step function as $\psi$ and piecewise
linear one as $\chi$, however we are going to perform these
approximations gradually. A step function is needed because we
only have a technical tool -- the Lemma \ref{L-main} --
established for this very case. 
A piecewise linear $\psi$ is not necessary, but
convenient. Eventually we will consider a
finite-dimensional ``discretized'' subset of
the set $\{\rho(X,\varphi)<\delta\}$ 
with appropriately chosen $\Delta$,  $X^\psi$, deterministic
curves $\psi, \chi$, and constants $\delta'_k$: in particular, we
will choose $\delta'_1 << \delta'_2 << \ldots <<
\delta'_{T/\Delta} << \delta$.
While performing all these approximations, we need to establish
simultaneously a special property: at any point $s$, the
Fenchel-Legendre adjoint to the $\dot \chi_s$ variable $\beta_s =
\beta_s[\psi_s,\dot \chi_s]$ (see below) can be chosen uniformly
bounded.

~

%\item %4
\noindent
{\bf 3.}
For any nonrandom curve $\psi \in C([0,T];E^d)$ -- although we will apply
this firstly to $\varphi$, but other functions are also
necessary for the analysis below -- we have, due to the
Lipschitz condition on
$f$,
\begin{equation}\label{eqC}
  \{\rho(X,\varphi)<\delta\} \supset
\{\rho(X^\psi,\chi)<\delta'\}
\end{equation}
if $\delta'$ and $\lambda : = \rho_{0,T}(\varphi,\psi)$ are small
enough with respect to $\delta$. (A small constant $\lambda>0$ is
used just within this step.) E.g., $ \delta' < \delta(e^{CT}CT +
1)^{-1}/2, \quad \lambda < \delta(e^{CT}CT + 1)^{-1}/2 $ suffice,
see below. Indeed,
$$
X_t = x + \int_0^t f(X_s,Y_s)ds, \quad X^\psi_t = x + \int_0^t
f(\psi_s,Y_s)ds,
$$
thence,
\begin{eqnarray*}
|X_t-X^\psi_t| \le \int_0^t |f(X_s,Y_s)ds - f(\psi_s,Y_s)|ds \le
C\int_0^t |X_s-\psi_s|ds
\\ \\
\le C\int_0^t |X_s-X^\psi_s|ds + C\int_0^t |X^\psi_s-\chi_s|ds
+C\int_0^t |\chi_s - \psi_s|ds;
\end{eqnarray*}
so on the set $\{\rho(X^\psi,\chi)<\delta'\}$,
\begin{eqnarray*}
|X_t-X^\psi_t| \le C\int_0^t |X_s-X^\psi_s|ds + C\delta' t
+C\lambda t,
\end{eqnarray*}
and, moreover, for every $\omega
\in \{\rho(X^\psi,\chi)<\delta'\}$ and $0\le t\le T$,
$$
\sup_{0\le t'\le t} |X_{t'}-X^\psi_{t'}|(\omega) \le
C\int_0^t \sup_{0\le
s'\le s} |X_{s'}-X^\psi_{s'}|(\omega)\,ds + C(\delta' +
\lambda) t.
$$
Since all SDE solutions $X_{t}$, $X^\psi_{t}$ are
continuous, $\sup_{0\le t'\le T}
|X_{t'}-X^\psi_{t'}| <\infty$  for each $\omega \in \Omega$.
By the standard ``non-random'' Gronwall inequality this implies 
that on the same set
$\{\rho(X^\psi,\chi)<\delta'\}$,
$$
\rho(X,X^\psi)(\omega) \le e^{CT}C(\delta'+\lambda)T.
$$
Now, still for any $\omega \in
\{\rho(X^\psi,\chi)<\delta'\}$,
\begin{eqnarray*}
\rho(X,\varphi)(\omega) \le \rho(X,X^\psi)(\omega) +
\rho(X^\psi,\chi)(\omega) + \rho(\chi,\varphi) 
 \\ \\
\le e^{CT}C(\delta'+\lambda)T + \delta' + \lambda =
(\delta'+\lambda)(e^{CT}CT + 1).
\end{eqnarray*}
Therefore, (\ref{eqC}) holds true. For example,
$$
\delta' < \delta(e^{CT}CT + 1)^{-1}/2, \quad \lambda <
\delta(e^{CT}CT + 1)^{-1}/2
$$
suffice. %{\it However, the smaller the better}.
In particular, it is true that
$$%\begin{equation}\label{eqC}
  \{\rho(X,\varphi)<\delta\} \supset
\{\rho(X^\varphi,\chi)<\delta'\},
$$%\end{equation}
if $\delta'$ and $\lambda$ are small enough with respect to
$\delta$. This bound will be used while establishing a lower
bound.

~

%\item %4
\noindent
{\bf 4.}
While establishing an upper bound, an opposite inclusion will be
useful,
\begin{equation}\label{eqCa}
  \{\rho(X,\varphi)<\delta\} \subset
  \{\rho(X^\psi,\chi)<2\delta(KT+1)\},
\end{equation}
if $\lambda :=
\max\left(\rho(\varphi,\psi),\rho(\varphi,\chi)\right) \le
\delta$. Indeed,
\begin{eqnarray*}
|X_t-X^\psi_t| \le \int_0^t |f(X_s,Y_s)ds - f(\psi_s,Y_s)|ds \le
K\int_0^t |X_s-\psi_s|ds
\\ \\
\le K\int_0^t |X_s-\varphi_s|ds + K\int_0^t |\psi_s-\varphi_s|ds;
\end{eqnarray*}
so on the set $\{\rho(X,\varphi)<\delta\}$,
\begin{eqnarray*}
|X_t-X^\psi_t| \le K\delta t +K\lambda t,
\end{eqnarray*}
and, moreover, on the same set,
$$
\rho(X,X^\psi) \le K(\delta+\lambda)T.
$$
Now, (\ref{eqCa}) follows from the inequalities,
\begin{eqnarray*}
&\rho(X^\psi,\chi) \le \rho(X,X^\psi) + \rho(X,\varphi) +
\rho(\chi,\varphi) \\ \\
&\le K(\delta+\lambda)T + \delta + \lambda.
\end{eqnarray*}

\noindent
{\bf 5.}
%\item 
Our next goal is the choice of appropriate functions $\chi$
and $\psi$. It is essential to keep the integral
$\int\limits_0^T L(\varphi_s,\dot\chi_s)\,ds$ close to
$S(\varphi)$. Also, by technical reasons we want some
discretization. Hence, we will use a trick well-known in
the definition of stochastic integrals based on the
following Lemma.

\begin{Lemma}\label{Le-int}
Suppose $g \in  L_1([0,T];R^d)$ and let
$\kappa_m(a):= [2^m a]2^{-m}$. Then there exists a
sequence $m'\to \infty$ such that for almost every $a\in
[0,1]$, 
\begin{equation}\label{lemmal1}
 \int_{0}^T |g(s) - g(\kappa_{m'}(s+a)-a)|\, ds \to 0,
\quad m'\to \infty.
\end{equation}
 
\end{Lemma}
For the proof for $g\in L_2([0,T];R^d)$ see \cite[Theorem
2.8.2]{Krylov (2002)}, however, for $L_1([0,T];R^d)$  the
proof
practically does not change: we approximate $g$ by
continuous functions $g_n$ -- which are dense in 
$L_1([0,T];R^d)$ -- 
and integrate with respect to $a\in [0,1]$. Then, for each
$g_n\in
C([0,T];R^d)$ the statement follows for every $a$ and for
the limiting function $g$ the assertion (\ref{lemmal1})
follows for almost every $\omega$ over some subsequence, 
as required.

Hence, applying this Lemma we may fix some $a\in
[0,1]$ for which there exists a sequence
$m'\to \infty$ 
such that 
%for almost every $a\in [0,1]$, 
\begin{equation}\label{L}
 \int_0^T |L(\varphi_{s}, \dot \varphi_{s}) -
L(\varphi_{\kappa_{m'}(s+a)-a}, \dot
\varphi_{\kappa_{m'}(s+a)-a})|\, ds \to 0,  \quad m'\to
\infty.
\end{equation}
Simultaneously for almost every $a\in [0,1]$, by virtue of
the same Lemma and because $\varphi$ is absolutely
continuous, we also have,
\begin{equation}\label{L2}
 \int_0^T |\dot \varphi_{s} - \dot
\varphi_{\kappa_{m'}(s+a)-a})|\, ds \to 0,  \quad m'\to
\infty,
\end{equation}
and 
\begin{equation}\label{L3}
 \sup_{0\le t \le T} |\varphi_{t} -
\varphi_{\kappa_{m'}(t+a)-a}| \to 0,  \quad m'\to
\infty, 
\end{equation}
each time over a new subsequence. Yet, to simplify
notations, in the sequel $m'$ will be replaced by $m$.
Denote 
$$
\psi_t^{m} = \psi_t := \varphi_{\kappa_{m}(t+a)-a}, \quad 
\dot \chi_t^{m} = \dot \chi_t :=
\dot\varphi_{\kappa_{m}(t+a)-a}, \quad \chi_t^{m} = \chi_t
:= \varphi_0 + \int_0^t \dot \chi_s\, ds.
$$
Notice that $\psi$ is piecewise constant (step function)
with finitely many values, while $\chi$ is piecewise
linear with finitely many values of slopes.

Let 
$$
{\cal S}^{\varphi}(\chi) := \int_0^T L(\varphi_s,
\dot\chi_s)\,ds. 
$$
Notice that ${\cal S}^{\varphi}(\varphi) = {\cal
S}(\varphi)$. 
Then (\ref{L}) implies
\begin{equation}\label{L4}
 |{\cal S}(\varphi) - {\cal S}^{\psi}(\chi)|\to 0,  \quad
m'\to \infty.
\end{equation}
At the same time we have, 
\begin{equation}\label{eqCb}
\fbox{ \mbox{$ \{\rho(X,\varphi)<\delta\} \supset
\{\rho(X,\psi)<\delta/2\}$,}}
\end{equation}
if $m$ is large enough. Moreover, in addition,
\begin{equation}%\label{eqCb} 
\fbox{ \mbox{$
\{\rho(X,\psi)<\delta/2\} \supset
\{\rho(X^{\varphi},\chi) < \tilde\delta'\}$,}}
\end{equation}
if $\tilde\delta'$ and $\lambda = \rho(\varphi,\chi)$
are
small enough with respect to $\delta$; hence, we can {\it
fix the
value $\tilde\delta'$ here.}

~

So, we can choose the functions $\psi$ and $\chi$ so that, 
firstly, $2^{-m} \le \Delta(\nu)$ (a value from
the Lemma \ref{L-main}); secondly, 
\begin{equation}%\label{s-s}
\left|S^\psi(\chi) - S(\varphi) \right| \le \nu;
\end{equation}
%thirdly, 
%\begin{equation}%\label{eqD}
%\fbox{ \mbox{$ \{\rho(X^{\varphi},\varphi^b
%)<\tilde\delta'\}
%\supset \{\rho(X^{\varphi}, \chi)<\tilde\delta'\times
%9/10\}$,}}
%\end{equation}
and, finally (see above (\ref{eqC})), if $\delta'$ is small
enough then also 
\begin{equation}%\label{eqE}
\fbox{ \mbox{$ \{\rho(X^{\varphi},\chi )<\tilde\delta'
%\times 9/10
\} \supset \{\rho(X^{\psi}, \chi)<\delta'\}$;}}
\end{equation}
for the latter we need only $\rho(\varphi,\chi)+\rho(\varphi,\psi)$ to
be small enough. 
%which means, in particular, $m$ large enough.

%Notice that the vector--function $\dot\chi$ has finitely
%many values and all of them are finite vectors satisfying
%$L(\psi_s, \dot\chi_s)<\infty$. Denote
%$$
%B:= \max |\dot\chi|, ??????????? (nado??)
%$$
%(NOT $b$) where $|\cdot|$ stands for the Euclidean norm. 

~

\noindent
{\bf 6.}
Suppose for some $s\in [0,T]$, the set $\{\alpha: \,
L(\psi_s,\alpha)<\infty\}$ has a non-empty interior 
with respect to its linear hull ${\cal L}[f,\psi_s]$,
that is, to the minial linear subspace containing $\{\alpha:
\, L(\psi_s,\alpha)<\infty\}$. For this interior --
non-empty or empty -- we will use notation ${\cal
L}^\circ[f,\psi_s]$. Since
$L(\psi_s,\dot\chi_s)<\infty$, this value is
attained as a
$\,\liminf\,$ of the values $L(\psi_s,\alpha)$, $\alpha
\in {\cal L}^\circ[f,\psi_s]$,
as $\alpha\to\dot\chi$ in the case ${\cal L}^\circ[f,\psi_s]
\not= \emptyset $, see 
%Freidlin and Wentzell and/or
\cite{Rockafellar (1970)}. It is a property of any such
$\alpha$ that
there exists a finite adjoint vector $\beta = \argmax_\beta
(\alpha\beta - H(\psi_s,\beta))$ given $\alpha$,
although this
adjoint may not be necessarily unique which we will discuss
shortly. Notice that, in particular, we have
$$
H(\psi_s,\beta) = (\alpha\beta - L(\psi_s,\alpha)),
\;\;
\mbox{as well as}\;\; L(\psi_s,\alpha) = (\alpha\beta -
H(\psi_s,\beta)).
$$
We can choose a vector $\dot{\tilde\chi_s}:=
\alpha \in {\cal L}^\circ[f,\psi_s]$ so that the value
$L(\psi_s,\dot{\tilde\chi_s})$ is close enough to
$L(\psi_s,\dot{\chi_s})$. Recall that there are finitely
many vector-values of $\dot{\chi_s}$ for any given $m$ and
$a$; correspondingly, we will choose finitely many
approximations satisfying $\dot{\tilde\chi_s}\in {\cal
L}^\circ[f,\psi_s]$. Let us also choose some
adjoint $\beta$ for each $\alpha = \dot{\tilde\chi_s}$ and
denote it by $\beta[\psi_s,\dot{\tilde\chi_s}]$. 

In the case if the set ${\cal L}^\circ[f,\varphi_s]$ is
empty, the function $H(\varphi_s, \beta)$ is linear in
$\beta$ and one can choose $\dot{\tilde\varphi_s} :=
\dot{\varphi_s}$ and 
$\beta[\varphi_s,\dot{\tilde\varphi_s}]=0$, see Appendix~A.

Notice that whatever is the case -- the interior ${\cal
L}^\circ[f,\psi_s]$ empty or
not -- and whatever is the choice of $\beta$ -- if not
unique -- in all cases there are finitely many of vectors
$\beta[\psi_s,\dot{\tilde\chi_s}]$ chosen. Hence, we
may denote 
\begin{equation}\label{eq-b}
\max_{0\le s\le T}|\beta[\psi_s,\dot{\tilde\chi_s}]|=: b < \infty.
\end{equation}
Notice that this value is fixed from now on. Let
$$
\tilde\chi_t := x + \int_0^t
\dot{\tilde\chi_s}\,ds, 
\quad S^{\psi}(\tilde\chi): = \int\limits_0^T
L(\psi_s,\dot{\tilde\chi_s})\,ds.
$$
We may assume that $\tilde\chi$ is as
close to $\varphi$ as we
like, say, $\rho(\tilde\chi, \psi) < \nu/3$ and also
\begin{equation}\label{ne1}
\left|S^{\psi}(\tilde\chi) - S(\varphi) \right| \le
\nu/3.
\end{equation}

~

%\item %6
\noindent
{\bf 7.} 
In the general case, the discretisations of $\varphi$ should
be read $\varphi^{\Delta, a}
= (\varphi_{\Delta-\tilde a},\varphi_{2\Delta
-\tilde a}, \ldots,\varphi_{m\Delta -\tilde a},
\varphi_{T})$, where $\tilde a = a -
[a/\Delta]\Delta$; if $a=0$ then we may use the
approximation $\varphi^\Delta = (\varphi_\Delta,
\varphi_{2\Delta}, \ldots, \varphi_{m\Delta})$,
$m\Delta=T$. Notice that '\emph{almost every} value' of
$a$ does not guarantee any particular value, so that we
cannot be sure about taking $a=0$. Hence, let us consider
the general case here. Denote ${k\Delta -\tilde a}=: t_k,
\; 1 \le k \le m,$ and $t_{m+1}:= T$ in the case of $\tilde
a \not = 0$ (and no $t_{m+1}$ in the case of $\tilde
a = 0$).

Since the drift of the diffusion $X^\psi$ is bounded -- 
$\|f\|_C < \infty$ -- we have straight away (however, cf. 
\cite[proof of the Lemma 7.5.1]{Freidlin and Wentzell
(1984)}),
\begin{equation}\label{eqGG}
\fbox{$\{\rho(X^\psi,\chi)<\delta'\} \supset
\{\rho((X^\psi)^{\Delta, a},\chi^{\Delta, a})<\delta''\},$}
\end{equation}
if $\delta''$ and $\Delta$ are small enough,
\begin{equation}\label{discr}
\delta'' < \delta''(\delta') \quad \mbox{and} \quad \Delta\le
\Delta(\delta')
\end{equation}
(notice that here $\Delta\le \Delta(\delta'')$ is \emph{not}
required), and assuming all our curves start at $x_0$ at
time
zero (hence, we do not include the starting point into the
definition of $\varphi^\Delta$). Here for discretized curves
we use the metric,
$$
\rho(\psi^{\Delta, a},\chi^{\Delta, a}) := \sup_{k}
|\psi_{t_k}-\chi_{t_k}|.
$$
Now, we are going to estimate from below the value in the right
hand side of the inequality,
\begin{equation}\label{delta''}%$$
\fbox{$P(\rho((X^\psi)^{\Delta, a},\chi^{\Delta, a}) <
\delta'') \ge E
\prod_{k} I( |X^\psi_{t_k} - \chi_{t_k}| <
\delta'_i),$}
\end{equation}%$$
where $\delta'_1<\delta'_2<\ldots < \delta'_{m+1} =
\min(\delta(\nu),\delta'')$, % =: \tilde\delta$,
$i = 1,\ldots, m$,
and $\delta(\nu)$ is from the Lemma~5; here all values $\delta'_i$
and certain auxiliary values $z_i$ will be chosen in the next two
steps as follows:
$$
m_{\nabla_\beta \tilde H}(\delta'_{k-1}+z_{k-1}) + \frac{\kappa}{2}\,
\delta'_{k-1}\le \frac{\kappa}{2}\,\delta'_k, \quad \& \quad
\delta'_{k-1}\le \frac{\delta'_k}{2}, \quad
 \& \quad m_{\tilde H}(\delta'_{k}) \le \nu,
$$
where $0<\kappa\le 1$. Emphasize that $\delta''$ and
$\Delta$ may be chosen arbitrarily small at this stage; in
particular, we require that they should satisfy the
conditions of the Lemma \ref{L-main}, which will be used
in the sequel, that is, we do require
$\delta''\le \delta(\nu)$ and $\Delta\le \Delta(\nu)$.
Hence,
{\it both $\delta''$ and $\Delta$ are fixed at this stage.}

~

%\item\label{item11} 
\noindent
{\bf 8.}\label{it11}
Now everything is prepared for the lower
estimate. We start with the estimation of the conditional
expectation $E(I(|X^\psi_{t_{m+1}} - \chi_{t_{m+1}}| <
\delta'_{m+1}) \mid  {\cal F}_{t_{m}})$ on the set
$\{|X^\psi_{t_{m}} -
\chi_{t_{m}}| < \delta'_{m} \}$. Let us apply the Cram\'er
transformation of measure. Let $|\beta|\leq b$, we will
choose
this vector a bit later (as
$\beta[\psi_{t_m},\dot\chi_{t_m+}]$). We get,
\begin{eqnarray*}
&E\left( I(|X^\psi_{t_{m+1}} - \chi_{t_{m+1}}| <
\delta'_{m+1})
| {\cal F}_{t_{m}} \right) = E^\beta \left(I(|X^\psi_{t_{m+1}} -
\chi_{t_{m+1}}| < \delta'_{m+1}) \times \right.
 \\ \\
&\left. \times \exp\left( -\eps^{-2} \beta (X^\psi_{t_{m+1}}
-X^\psi_{t_{m}}) +
\eps^{-2}\Delta_m H^{\eps,\psi}_m(X^\psi_{t_{m}},\beta)
\right)| {\cal F}_{t_{m}} \right),
\end{eqnarray*}
where $E^\beta$ is the (conditional) expectation with respect to
the measure $P^\beta$ defined on the sigma-field
$ {\cal F}_{t_{m+1}}$ given $ {\cal F}_{t_{m}}$, by its density
$$
  \frac{dP^\beta}{dP}(\omega) =
\exp\left(\eps^{-2} \beta (X^\psi_{t_{m+1}} 
- X^\psi_{t_{m}}) - \eps^{-2}\Delta_m
\tilde H^{\eps,\psi}_m (X^\psi_{t_{m}},\beta)\right),
$$
where $\Delta_m = t_{m+1} - t_{m}$ (and later on, $\Delta_k
= t_{k+1} - t_{k}$; notice that all $\Delta_k \le \Delta$)
$$
\eps^{-2}\Delta_m \tilde 
H^{\eps,\psi}_m(X^\psi_{t_{m}},\beta) := \log E
\left(\exp\left( \eps^{-2}\beta (X^\psi_{t_{m+1}}
-X^\psi_{t_{m}})\right) | {\cal F}_{t_m} \right).
$$
Notce that by virtue of the Lemma \ref{L-main}, 
$$
\tilde H^{\eps,\psi}_m(X^\psi_{t_{m}},\beta)
\to \tilde H(\psi_{t_{m}}, X^\psi_{t_{m}},\beta), 
\quad \eps \to 0, 
$$
uniformly over $|\beta| \le b$. Indeed, by definition of 
$X^\psi$, 
$$
X^\psi_{t_{m+1}} - X^\psi_{t_{m}} = 
\int^{t_{m+1}}_{t_{m}} f(\psi_{t_{m}}, Y_s)\,ds.
$$
Thus, the inequality (\ref{eqDD}) of the Lemma \ref{L-main} implies, 
\begin{eqnarray*}
 |\tilde H^{\eps,\psi}_m(X^\psi_{t_{m}},\beta)
- \tilde H(\psi_{t_{m}}, X^\psi_{t_{m}},\beta)|
\hspace{3cm}
 \\\\
= \left|
\eps^{2}\Delta_m^{-1} \log E
\left(\exp\left( \eps^{-2}\beta (X^\psi_{t_{m+1}}
-X^\psi_{t_{m}})\right) | {\cal F}_{t_m}\right) - 
\tilde H(\psi_{t_{m}}, X^\psi_{t_{m}},\beta)
\right|\hspace{1cm}
 \\\\
= \left|
\eps^{2}\Delta_m^{-1} \log E
\left(\exp\left( \eps^{-2}\beta 
\int^{t_{m}+\Delta_m}_{t_{m}} f(\psi_{t_{m}}, Y_s)\,ds\right) 
\mid {\cal F}_{t_m}\right) - 
\tilde H(\psi_{t_{m}}, X^\psi_{t_{m}},\beta)
\right| 
 \\\\
\le  \nu \eps^{-2}\Delta_m^{}. \hspace{5.5cm}
\end{eqnarray*}
Also notice that on the set
$\{|X^\psi_{t_m} - \chi_{t_m}| < \delta'_{m}\}$ we have
$$
\eps^{-2}\left|\beta(X^\psi_{t_{m}} 
- \chi_{t_{m}})\right| 
\le \eps^{-2} \, b\, \delta'_{m}
$$
and on the set the set
$\{|X^\psi_{t_{m+1}} -\chi_{t_{m+1}}| <
\delta'_{m+1}\}$,
$$
\eps^{-2}\left|\beta(X^\psi_{t_{m+1}} 
- \chi_{t_{m+1}})\right| 
\le \eps^{-2} \, b\, \delta'_{m+1}.
$$
Hence, for $\eps>0$ small enough on the set
$\{|X^\psi_{t_m} - \chi_{t_m}| < \delta'_{m} \}$
we estimate,
\begin{eqnarray}\label{eq17}
& E\left[ I(|X^\psi_{t_{m+1}} -\chi_{t_{m+1}}| <
\delta'_{m+1}) \mid  {\cal F}_{t_{m}}\right] \nonumber
 \\ \nonumber \\
&=E^\beta \left[ I(|X^\psi_{t_{m+1}} -\chi_{t_{m+1}}| <
\delta'_{m+1}) \right. 
 \nonumber \\ \left.\right.\nonumber \\ 
&\left. \times
\exp\left(\eps^{-2}\beta(X^\psi_{t_{m+1}}
 - X^\psi_{t_{m}}) - \eps^{-2}\Delta_m
\tilde H^{\eps,\psi}_m(X^\psi_{t_{m}},\beta)\right) \mid
 {\cal F}_{t_{m}} \right]
 \nonumber \\ \nonumber \\
& \geq E^\beta \left( I(|X^\psi_{t_{m+1}} -\chi_{t_{m+1}}| <
\delta'_{m+1}) \exp\left( -\eps^{-2} \Delta_m \beta
\left((\chi_{t_{m+1}}-\chi_{t_{m}})/\Delta_m \right)
 \nonumber \right. \right. \\ \nonumber \\
& \left. \left. - \frac{\Delta_m}{\eps^{2}} (\tilde H(\psi_{t_{m}},
X^\psi_{t_{m}},\beta) + \nu) -
\frac{b\,(\delta'_{m+1}+\delta'_{m})}{\eps^{2}} \right) |
 {\cal F}_{t_{m}}\right).
\end{eqnarray}

~

\noindent
Now, let us choose $\beta = \beta(m+1) =
\beta[\psi_{t_{m}},\dot\chi_{t_{m}+}] \; =
\argmax_{\beta}(\beta\dot\chi_{t_{m}+} -
H(\psi_{t_{m}},\beta))$. As was explained above,
$|\beta(m+1)|\leq b$ and, moreover,
\begin{eqnarray*}
\beta(m+1)\dot\chi_{t_{m}+} -
H(\psi_{t_{m}},\beta(m+1)) %\\\\
= L(\psi_{t_{m}}, \dot\chi_{t_{m}+}),
\end{eqnarray*}
and
\begin{equation}\label{chiprime}
\dot\chi_{t_{m}+} = \nabla_\beta
H(\psi_{t_{m}},\beta(m+1)).
\end{equation}
So (\ref{eq17}) implies (with $\beta=\beta(m+1)$),
\begin{eqnarray}\label{eq-low}
&E\left( I(|X^\psi_{t_{m+1}} -\chi_{t_{m+1}}| <
\delta'_{m+1}) \mid {\cal F}_{t_{m}}\right)
 \nonumber \\\nonumber \\
&\geq \exp\left(-\eps^{-2}\Delta_m (L(\psi_{t_{m}},
\dot\chi_{t_{m}+})+\nu) -
b \eps^{-2}(\delta'_{m+1} + \delta'_{m})
\right)\times
 \nonumber \\\nonumber \\
&\times  \exp\left(- \eps^{-2}\Delta_m (\tilde H(\psi_{t_{m}},
X^\psi_{t_{m}},\beta) -
\tilde H(\psi_{t_{m}},\psi_{t_{m}},\beta))\right)
 \nonumber  \\ \nonumber \\
&\times E^{\beta(m+1)}\left( I(|X^\psi_{t_{m+1}} -
\chi_{t_{m+1}}| < \delta'_{m+1}) \mid  {\cal F}_{t_{m}}\right)
 \nonumber \\\nonumber \\
&\geq
\exp\left(-\eps^{-2}\Delta_m \left(L(\psi_{t_{m}},
\dot\chi_{t_{m}+}) + 2\nu\right)  - b \eps^{-2}
(\delta'_{m+1} + \delta'_{m}) \right)\times
 \nonumber  \\ \nonumber \\
&\times E^{\beta(m+1)}\left( I(|X^\psi_{t_{m+1}} -
\chi_{t_{m+1}}| < \delta'_{m+1})| {\cal F}_{t_{m}}\right).
\end{eqnarray}
We have used uniform continuity of $\tilde H(x,\cdot,\beta)$
over $|\beta|\le b$ and $x~\in~R^d$:
\begin{eqnarray*}
|\tilde H(\psi_{t_{m}}, X^\psi_{t_{m}},\beta) -
\tilde H(\psi_{t_{m}},\psi_{t_{m}},\beta)| \\\\
\le m_{\tilde H}(|X^\psi_{t_{m}} - \psi_{t_{m}}|)\le
m_{\tilde H}(\delta'_{m}) \le \nu
\end{eqnarray*}
on the set $|X^\psi_{t_{m}} - \psi_{t_{m}}|\le 
\delta'_{m}$ 
(recall that here $m_{\tilde H}$ stands for the modulus of continuity
of $\tilde H$ for $|\beta|\le b$), as $\delta'_{m}$ is small
enough.

~

%\item\label{item12} 
\noindent
{\bf 9.}\label{it12} 
Let us show that given $\delta'_{m+1}$, there exists 
$C_{m+1}>0$ such that on the set $\{|X^\psi_{t_{m}} -
\chi_{t_{m}}| < \delta'_{m}\}$, 
\begin{equation}\label{eqLLN}
  E^{\beta(m+1)}\left( I(|X^\psi_{t_{m+1}} - \chi_{t_{m+1}}|
< \delta'_{m+1}) \mid  {\cal F}_{t_{m}}\right) \geq 1 -
\exp(-C_{m+1} \eps^{-2}),
\end{equation}
if $\eps$ is small enough.
There exists a finite number of vectors $v_1,\, v_2,\,
\ldots, v_{2d}$ such that $\|v_k\|=1
\; \forall k$ (any orthonormal basis would do accomplished
by its ``symmetric'' transformation, i.e. with each
coordinate vector $v$ we consider $-v$ as well), and for any 
(non-random) vector $\xi$ and any positive $c$, 
$$
|\xi| > c \quad \Longrightarrow \quad \exists\;\; 1\le k\le 2d:
\quad \xi v_k > \kappa c,
$$
where $\kappa = (1/d)^{1/2}$ (notice that $\kappa\le 1$).
Then, 
\begin{eqnarray*}
& E^{\beta(m+1)} (I(|X^\psi_{t_{m+1}} - \chi_{t_{m+1}}| >
\delta'_m) \mid {\cal F}_{t_{m}})
 \\ \\
& \leq \sum^{2d}_{k=1} E^{\beta(m+1)} (I((X^\psi_{t_{m+1}}
- X^\psi_{t_{m}} - \chi_{t_{m+1}} + \chi_{t_{m}}) v_k
 \\ \\
& > \kappa (\delta'_{m+1} - \delta'_{m}) \mid
{\cal F}_{t_{m}}),
\end{eqnarray*}
given $\{|X^\psi_{t_{m}} - \chi_{t_{m}}| <
\delta'_{m}\}$. 
Let $\nu\,'_{m}>0$ (this is a new constant which has nothing
to do with $\nu$ and will be fixed shortly, see (\ref{nure}) 
below; we need it only
while establishing the inequality (\ref{eqLLN})). 
By exponential Chebyshev's inequality we estimate, 
for any $v := v_k$ and any $0\leq z\leq 1$
%(actually, any $z\ge 0$ would do, but we need to restrict
%the modulus of the variable $\beta + vz$ by $b+1$) 
on the set $\{|X^\psi_{t_{m}} - \chi_{t_{m}}| 
< \delta'_{m}\}$,
\begin{eqnarray}\label{eq20}
&E^{\beta(m+1)}\left(
I((X^\psi_{t_{m+1}} - X^\psi_{t_{m}}
- \chi_{t_{m+1}} + \chi_{t_{m}}) v
> \kappa
(\delta'_{m+1} - \delta'_{m}))|{\cal F}_{t_{m}}\right)
 \nonumber \\ \nonumber \\
& = E^{\beta(m+1)}\left(
I(z \eps^{-2}\, (X^\psi_{t_{m+1}} - X^\psi_{t_{m}}
- \chi_{t_{m+1}} + \chi_{t_{m}}) v
  \right. \nonumber \\ \nonumber \\ & \left. 
 >  z \eps^{-2}\, \kappa
(\delta'_{m+1} - \delta'_{m}))|{\cal F}_{t_{m}}
\phantom{X^\psi_{t}} \!\!\!\!\!\!\!\!\!\right)
 \nonumber \\ \nonumber \\
&\leq \exp(-(\delta'_{m+1} 
 - \delta'_{m})z \kappa \eps^{-2})
  \nonumber \\ \nonumber \\
& \times E^{\beta(m+1)} 
\exp( z \eps^{-2}\, (X^\psi_{t_{m+1}} - X^\psi_{t_{m}}
- \chi_{t_{m+1}} + \chi_{t_{m}}) v)
\nonumber \\ \nonumber \\
&\leq \exp(-(\delta'_{m+1} 
 - \delta'_{m}) z \kappa
\eps^{-2})\exp\left(\eps^{-2} [- z
v \dot\chi_{t_{m}+} \Delta_m  \right.
 \nonumber \\ \nonumber \\
&\left. + \tilde H^{\eps,\psi}(X^\psi_{t_{m}},
\beta(m+1)+ v z) -
\tilde H^{\eps,\psi}(X^\psi_{t_{m}},\beta(m+1))
+ 2\nu\,'_{m-1}]\right)
 \nonumber \\ \nonumber \\
&\leq \exp(-(\delta'_{m+1} - \delta'_{m})z \kappa
\eps^{-2})\exp\left(\eps^{-2} [- z
v\dot\chi_{t_{m}+} \Delta_m \right.
 \nonumber \\ \nonumber \\
&\left. + \tilde H(\psi_{t_{m}},X^\psi_{t_{m}},\beta(m+1) +
v z) \right.
 \nonumber \\ \nonumber \\
&\left. - \tilde H(\psi_{t_{m}},X^\psi_{t_{m}},\beta(m+1))
+ 2\nu'_{m}]\right),
\end{eqnarray}
if $\eps$ is small enough. 
We used here the identity 
$\chi_{t_{m+1}} - \chi_{t_{m}} = \Delta_m \dot\chi_{t_{m}+}$.
Denote
\begin{eqnarray*}
&h(z) :=  (\delta'_{m+1} - \delta'_{m})\kappa z +
\dot\chi_{t_{m}+} v z \Delta_m
 \\ \\
&-[\tilde H(\psi_{t_{m}},X^\psi_{t_{m}},\beta(m+1)+ v z) -
\tilde H(\psi_{t_{m}},X^\psi_{t_{m}},\beta(m+1))],
\end{eqnarray*}
so that the rightmost side of (\ref{eq20}) may be 
represented as 
$$
\exp(- \eps^{-2} h(z)).
$$
Notice that $h(0)=0$. Moreover, since
$\dot\chi_{t_{m}+} = \nabla_\beta
\tilde H(\psi_{t_{m}},\psi_{t_{m}},\beta(m+1))$ 
(see (\ref{chiprime})), we have on the set 
$\{|X^\psi_{t_{m}} - \chi_{t_{m}}| < \delta'_{m}\}$, 
\begin{eqnarray*}
& h'(0) = (\delta'_{m+1}-\delta'_{m}) \kappa +
\dot\chi_{t_{m}+}v \Delta_m - \nabla_\beta
\tilde H(\psi_{t_{m}},X^\psi_{t_{m}},\beta(m+1))v \Delta_m
 \\ \\
& = (\delta'_{m+1} - \delta'_{m})\kappa \Delta_m 
+ \nabla_\beta \tilde H(\psi_{t_{m}},\psi_{t_{m}},\beta(m+1))
v \Delta_m
 \\ \\
& - \nabla_\beta
\tilde H(\psi_{t_{m}},X^\psi_{t_{m}},\beta(m+1))v\Delta_m
 \\\\
& \ge (\delta'_{m+1} - \delta'_{m})\kappa -
m_{\nabla_{\beta} \tilde H}(\delta'_{m})\Delta =:C'_{m+1} > 0
\end{eqnarray*}
(recall that $\Delta_m \le \Delta$ and 
that here $m_{\nabla_{\beta} \tilde H}$ stands for 
the modulus of continuity of the
function $\nabla_\beta \tilde H$ given $|\beta(m)|\le b+1$ ($b+1$
will be useful in the sequel, although here $b$ would be
enough)). The inequality 
$C'_{m+1}=(\delta'_{m+1}-\delta'_{m})\kappa - m_{\nabla_{\beta}
\tilde H}(\delta'_{m})\Delta > 0$ holds true provided 
$\delta'_{m}$ is
small enough in comparison to $(\delta'_{m+1}-\delta'_{m})$,
e.g.,
$$
m_{\nabla_{\beta} \tilde H}(\delta'_{m})\Delta \le
\frac{\kappa}{2}(\delta'_{m+1}-\delta'_{m}),
$$
or, equivalently,
\begin{equation}\label{modulus1}
m_{\nabla_\beta \tilde H}(\delta'_{m}) \Delta + \frac{\kappa}{2}
\,\delta'_{m}\le \frac{\kappa}{2}\,\delta'_{m+1}.
\end{equation}
Recall that a slightly stronger assumption was used in the rule 
of choosing $\delta'_m$ and we will need a stronger version 
in a minute, see (\ref{modulus2}) below.

Moreover, since $\nabla_\beta \tilde H$ is bounded and continuous
due to the Lemma \ref{L-dif}, then $h'(z) \geq C_{m}/2$ for
small $z$, say, for $0\leq z\leq z_{m}$ (thus, $z_{m}$ is
fixed here), on the set $\{|X^\psi_{t_{m}} -
\chi_{t_{m}}| < \delta'_{m}\}$. Indeed,
\begin{eqnarray*}
h'(z) = (\delta'_{m+1}-\delta'_{m})\kappa  +
\dot\chi_{t_{m}+}v \Delta_m
 \\ \\
- \nabla_\beta
\tilde H(\psi_{t_{m}},X^\psi_{t_{m}},\beta(m+1)+vz)v \Delta_m
 \\ \\
= (\delta'_{m+1}-\delta'_{m})\kappa + \nabla_\beta
\tilde H(\psi_{t_{m}},\psi_{t_{m}},\beta(m+1))v\Delta_m
\\ \\
- \nabla_\beta
\tilde H(\psi_{t_{m}},X^\psi_{t_{m}},\beta(m+1)+vz)v\Delta_m
\\\\
\ge (\delta'_{m+1}-\delta'_{m})\kappa - m_{\nabla
\tilde H}(\delta'_{m}+z)\Delta. 
\end{eqnarray*}
So, $h(z_{m}) \geq C'_{m+1} z_{m} /2$, provided $z_{m}$
along with $\delta'_{m}$ are both small in comparison to
$(\delta'_{m+1}-\delta'_{m})$, for example, if
\begin{equation}\label{modulus2} 
m_{\nabla \tilde H}(\delta'_{m}+z_{m})\Delta \le
(\delta'_{m+1}-\delta'_{m})\kappa/2,
\end{equation}
rather than (\ref{modulus1}). Hence, under the assumption of
(\ref{modulus2}), the right hand side in (\ref{eq20})
with $z=z_m$ on the set $\{|X^\psi_{t_{m}} -
\chi_{t_{m}}| < \delta'_{m}\}$ does not                
exceed the value
$$
\exp(\eps^{-2}(2\nu\,'_{m}- h(z_m))) \leq \exp(-
C'_{m+1} z_{m}\eps^{-2}/4])
$$
if we choose
\begin{equation}\label{nure}
\nu\,'_{m} < C'_{m+1} z_{m}/8.
\end{equation}
Recall that the constant $\nu\,'_{m}$ should have been fixed in
the beginning of this step of the proof; hence, we can do it
now, once we have chosen $z_{m}$, since the latter does
not require any knowledge of $\nu\,'_{m}$. Given
$\{|X^\psi_{t_{m}} - \chi_{t_{m}}| < \delta'_{m}\}$, this
implies the bound, 
$$
E^{\beta(m+1)}\left( I(|X^\psi_{t_{m+1}} - \chi_{t_{m+1}}|
\geq
\delta'_{m+1})|F_{t_{m}}\right) \leq \exp(- C'_{m+1}
z_m \eps^{-2}/4),
$$
which is equivalent to (\ref{eqLLN}) with 
$C_{m+1}:= C'_{m+1} z_m$. 
In turn, (\ref{eqLLN}) implies the estimate
\begin{eqnarray*}
&P(|X^\psi_{t_{m+1}} - \chi_{t_{m+1}}| < \delta'_{m+1}|
 {\cal F}_{t_{m}})
 \nonumber \\ \nonumber \\
&\geq \exp\left(-\eps^{-2}\Delta_m(L(\psi_{t_m},
\dot\chi_{t_m+}) + 3\nu) -
b \eps^{-2}(\delta'_{m+1}+\delta'_{m}) \right),
\end{eqnarray*}
still on $\{|X^\psi_{t_{m}} - \chi_{t_{m}}| <
\delta'_{m}\}$, if $\eps$ is small enough. Indeed, $\nu$, 
$C_{m+1}$ and $\Delta_m$ being fixed, one can choose $\eps$ 
so that
$$
1-\exp(- C_{m+1} \eps^{-2}) \ge \exp(-1) \ge
\exp(-\nu(\Delta_m\eps^{-2})).
$$

~

%\item %13
\noindent
{\bf 10.}
By ``backward'' induction from $k=m$ to $k=1$,
choosing at
each step $\delta'_{k-1}$ and $z_{k-1}$ small enough in
comparison to
$\delta'_{k}-\delta'_{k-1}$,
\begin{equation}\label{modulus3}
m_{\nabla_\beta \tilde H}(\delta'_{k-1}+z_{k-1})\Delta + \frac{\kappa}{2}
\delta'_{k-1}\le \frac{\kappa}{2}\delta'_k, \;\; \& \;\;
\delta'_{k-1}\le\delta'_k/2, \;\; \& \;\; m_{\tilde H}(\delta'_{k-1})<\nu
\end{equation}
(cf. (\ref{modulus2})), as well as all auxiliary values $C_{k-1}$,
for $\varepsilon$ small enough and since 
$\sum_{k=1}^{m+1}\delta'_k \le 2 \delta'_{m+1}$, 
we get the desired lower bound: 
\begin{eqnarray*}
&P(| X^{\psi}_{t_{m+1}}-\varphi_{t_{m+1}}|
<\delta'_{m+1}, \ldots,
|X^{\psi}_{t_{1}}-\varphi_{t_{1}}|<\delta'_1)
 \\ \\
& \ge \exp\left(-\eps^{-2} \sum^m_{i=0}
(L(\psi_{(m-i)\Delta}, \dot\chi_{(m-i)\Delta+})+3\nu )\Delta_i -
2b\eps^{-2}\sum_{k=1}^{m+1}\delta'_k\right)
 \\ \\
&\ge \exp\left(-\eps^{-2}(\int_0^T L(\psi_s,\dot\chi_s)\,ds
+ 3\nu T)
- 4b\eps^{-2}\delta'_{m+1}\right)
 \\\\
&\ge \exp \left(-\eps^{-2}(S_{0T}(\varphi ) + \nu
(3T+2))\right),\qquad\eps \to 0,
\end{eqnarray*}
provided $4b\delta'_{m+1}<\nu$.
This is equivalent to (\ref{eq5a}). This bound is uniform in
$x\in E^d,\,
|y|\leq r$, and $\varphi\in\Phi_x(s)$ for any $r,s>0$,
similar to
the Lemma 7.4.1 from \cite{Freidlin and Wentzell (1984)}.

~

%\item %14
\noindent
{\bf 11.}
The property of the rate function $S$ to be a ``good rate
function'' can be shown as in \cite{Freidlin and Wentzell
(1984)}, using
the semi-continuity of the function $L(x,y)$ with respect
to $y$ and
continuity with respect to $x$ variable (see 
\cite[Lemma 7.4.2]{Freidlin and Wentzell (1984)}).

~

%\item %15
\noindent
{\bf 12.}
\underline{\bf Second part of the proof:}  the
upper bound. Assume that the assertion (\ref{eq4a}) is not
true, that is, there
exist $s$ and $\nu>0$ with the following properties:
$$
\forall \bar\delta>0, \;
\mbox{there exists} \;\delta_0<\bar\delta, \; \forall
\bar\varepsilon, \; \mbox{there exists} \;
\varepsilon<\bar\varepsilon:
$$
$$
P(\rho(X,\Phi_x(s))>\delta_0) > \exp(-\varepsilon^{-2}(s - \nu)).
$$
In other words, for some (hence, actually, for any)
$\delta_0>0$ arbitrarily close to zero, there exists a sequence
$\varepsilon_n\to 0$ such that
\begin{equation}\label{contra}%$$
P(\rho(X,\Phi_x(s))>\delta_0)
> \exp(-\varepsilon_n^{-2}(s - \nu)).
\end{equation}%$$
We fix any such $\delta_0>0$.

~

%\item %16
\noindent
{\bf 13.}
Since $f$ is bounded, all possible trajectories of $X^\psi$
for any 
$\psi$ belong to some compact $F\subset C[0,T;R^d]$. Due to
semi-continuity of the functional $S^\psi(\varphi)$ with
respect to  $\psi$,
for any $\nu>0$ there exists a value $\delta > 0$ such that
$\rho(\varphi,\psi)<\delta$ and
$S(\varphi)>s$ imply $S^\psi(\varphi)>s-\nu/2$. 
Hence, let us define for each $\varphi \in C[0,T; R^d]$ a
positive value (notice that this definition differs slightly
from that given in \cite{Freidlin and Wentzell (1984)}; for
the latter -- without $\sup$ -- there is no reason to be
necessarily semi-continuous)
$$
\delta_\nu(\varphi):= \sup(\delta: \,
\rho(\varphi,\psi)<\delta \; \mbox{and} \; 
S(\varphi)>s \; \Longrightarrow \;
S^\psi(\varphi)>s-\nu/2).
$$
Since $S^\psi(\varphi)$ is lower semi-continuous with
respect to $\varphi\,$, too, similarly to $S(\varphi)$,
then it follows that $\delta_\nu(\varphi)$ is also lower
semi-continuous with respect to $\varphi$. Indeed, 
let $\varphi^n \to \varphi$, $n\to\infty$; we ought to
show that $\liminf_{n\to\infty} \delta_\nu(\varphi^n) \ge
\delta_\nu(\varphi)$. 
We have, 
$$
\delta_\nu(\varphi^n):= \sup(\delta: \,
\rho(\varphi^n,\psi)<\delta \; \mbox{and} \; 
S(\varphi^n)>s \; \Longrightarrow \;
S^\psi(\varphi^n)>s-\nu/2).
$$
Suppose $0 < \bar \delta < \delta_\nu(\varphi)$ and
$\rho(\varphi^n,\psi) < \bar \delta$. We want to show
that $S^\psi(\varphi^n)>s-\nu/2$. 
Since $\rho(\varphi^n,\varphi)<\delta_\nu(\varphi) -
\bar \delta$ for $n$ large enough, then we also have
$\rho(\varphi,\psi)<\delta_\nu(\varphi)$. Then, by
definition of $\delta_\nu(\varphi)$,
$S^\psi(\varphi)>s-\nu/2$. Since by Fatou's lemma, 
$\liminf_{n\to\infty} S^\psi(\varphi^n)>s-\nu/2$, this
implies $S^\psi(\varphi^n)>s-\nu/2$ for $n$ large enough.
The latter signifies that, indeed,
$\liminf_{n\to\infty} \delta_\nu(\varphi^n) \ge
\delta_\nu(\varphi)$, that is, that $\delta_\nu(\varphi)$ is
lower-semicontinuous, as required.

Thus, as every lower
semi-continuous function, $\delta_\nu(\varphi)$ attains
its minimum on any compact and, hence, the minimum
over any compact must be positive. 

Further, consider $F_1$, the compact obtained from $F$ by
dropping the
$\delta_0/2$-neighbourhood of the set $\Phi_x(s) = \{\varphi
\in
C[0,T;R^d]: \, \varphi_0=x, \, S(\varphi)\le s\}$. Denote
$\bar\delta_\nu = \inf_{\varphi\in F_1} \delta_\nu(\varphi)$, and
take any $\delta' \le
\min\left(\bar\delta_\nu/(4KT+2),\delta_0/2\right)$ where $K$ is a
Lipschitz constant of $f$.
Choose a finite $\delta'$-net for the set $F_1$, let
$\varphi^1, \ldots,
\varphi^N$ be its elements. All of them do not belong to
$\Phi_x(s)$, hence, $S(\varphi^i)\ge s'$ with some $s'>s$.
Notice that
$$
\{\rho(X,\Phi_x(s))>\delta_0\} \subset \bigcup_{i=1}^N
\{\rho(X,\varphi^i)<\delta'\}.
$$
Then, 
for any $n$ there exists an index
$i$ such that
\begin{equation}\label{5.14}
P(\rho(X,\varphi^i)\le \delta') >
N^{-1}\exp(-\varepsilon_n^{-2}(s
- \nu)).
\end{equation}

There is a finite number of $i = 1, \ldots, N$. Thus, there
exists at least one $i$ such that (\ref{5.14}) holds true
for this $i$ for some subsequence $n'\to\infty$ and
correspondingly $\varepsilon_{n'}\to 0$; however, we will
keep the notation $n$ for simplicity.
We may rewrite (\ref{5.14}) as
\begin{equation}\label{5.14a}
P(\rho(X,\varphi^i)\le \delta') > \exp(-\varepsilon_n^{-2}(s -
\nu)),
\end{equation}
since $N$ does not depend on $\varepsilon_n$, strictly speaking
with some new $\nu>0$; however, it is again convenient to
keep the same
notation. Denote $\varphi(\delta'):= \varphi^i$ with
this~$i$ (any one if not unique).

~

%\item %18
\noindent
{\bf 14.}
Consider a \emph{sequence} $\delta'\to 0$ such that a
corresponding
function $\varphi(\delta')$ does exist for any $\delta'$ from this
sequence. Recall that $\delta_0$ is fixed. All these functions
satisfy inequality
$$
S(\varphi(\delta'))\ge s'> s,
$$
since $\rho(\varphi^i,\Phi_x(s))\ge \delta_0/2$. Also we
have, $S(\varphi(\delta'))<\infty$, which implies
$$
\sup_t |\dot\varphi_t(\delta')|\le C,
$$
because, due to the boundedness of $f$, function $L(x,
\alpha)$ equals infinity for every $|\alpha| > \|f\|_C$.
By virtue of the Arcela-Ascoli Theorem, it is
possible to extract from
this set of functions  a subsequence which converges in
$C[0,T;R^d]$ to some limit, $\bar\varphi$. Since
$\rho(\varphi(\delta'),\Phi_x(s))\ge \delta_0/2$, we have,
$\rho(\bar\varphi,\Phi_x(s))\ge \delta_0/2$, hence,
$$
S(\bar\varphi) > s,
$$
and, in particular, the lower bound (\ref{eq5a}) can be applied.
However, due to the construction, the function $\bar\varphi$
satisfies one more lower bound,
\begin{equation}\label{better}%$$
\liminf_{\delta'\to 0} \limsup_{\varepsilon\to 0} \varepsilon^2\ln
P(\rho(X,\bar\varphi)<\delta') \ge - s+\nu.
\end{equation}%$$
Indeed, the latter follows from (\ref{5.14a}) because, e.g.,
$$
P\left(\rho(X,\bar\varphi)\le
\delta'+\rho(\bar\varphi,\varphi(\delta'))\right) \ge
P(\rho(X,\varphi(\delta'))\le \delta') >
\exp(-\varepsilon_n^{-2}(s - \nu)).
$$

Due to (\ref{better}), there exists $\hat\delta'>0$ such that for
smaller $\delta'$'s (a sequence)
$$
\limsup_{\varepsilon\to 0} \varepsilon^2\ln P(\rho(X,\bar
\varphi)<\delta') \ge - s+\nu/2.
$$
In fact, this implies the same inequality for any $\delta'>0$,
because with any $\delta'$ for which the inequality holds true,
every greater value would do as well. Therefore, for any
$\delta'$,
there exists $\varepsilon>0$ (arbitrarily small) such that

\begin{equation}\label{ge}%$$
\varepsilon^2\ln P(\rho(X,\bar\varphi)<\delta') \ge - s +\nu/3 =
-(s -\nu/3).
\end{equation}%$$
We are going to show that this leads to a contradiction.

~

%\item %?
\noindent
{\bf 15.}
Consider the case $S(\bar\varphi)<\infty$. Remind that
$ S(\bar\varphi) > s. $ Denote
\begin{eqnarray*}%$$
L^b(x,y) = \sup_{|\beta|\le b}(\beta y - H(x,\beta)),
\\\\
\ell^b(x,y) := L(x,y) - L^b(x,y)
\\\\
\equiv \sup_{\beta}(\beta y - H(x,\beta)) - \sup_{|\beta|\le
b}(\beta y - H(x,\beta)).
\end{eqnarray*}
Consider the function $\ell^b(\bar\varphi_t,\dot{\bar\varphi}_t)$.
We have, $$0\le \ell^b(\bar\varphi_t,\dot{\bar\varphi}_t)\le
L(\bar\varphi_t,\dot{\bar\varphi}_t).
$$
Moreover,
$$
\ell^b(\bar\varphi_t,\dot{\bar\varphi}_t)\to 0, \quad b\to\infty,
$$
and the function $\ell$ is decreasing with $b\to\infty$. Hence,
given $\nu>0$, one can choose a $b>0$ such that
$$
\int_0^T \ell^b(\bar\varphi_t,\dot{\bar\varphi}_t)\,dt < \nu/20.
$$
Notice that we have chosen $b$, which is now fixed for the second 
part of the proof of the Theorem. Moreover, one can also choose a
discretisation step $\Delta$ (see above, step 5 of
the proof and, in particular, the Lemma \ref{Le-int}) such
that for almost every $a \in [0,1]$
$$
\int_0^T
\ell^b(\bar\varphi_{\kappa_m(t+a)-a},\dot{\bar\varphi}_{\kappa_m(t+a)-a})\,dt
< \nu/10,
$$
and, correspondingly,
\begin{equation}\label{elb}
\int_0^T
L^b(\bar\varphi_{\kappa_m(t+a)-a},\dot{\bar\varphi}_{\kappa_m(t+a)-a})\,dt
> s - \nu/10.
\end{equation}
%again assume for simplicity of presentation that $a=0$. 
In addition, we require $\Delta\le \Delta(\nu/20)$ (see the
Lemma \ref{L-main}). Hence, we have chosen $\Delta$ and
$m=T/\Delta$. We also fix any $a\in [0,1]$ satisfying
(\ref{elb}).

~

%\item %19
\noindent
{\bf 16.}
Further, let 
%with $a=0$ for simplicity, let
$$
%\psi_t := \bar\varphi_{\kappa_m(t)}, \quad 
%\dot\chi_t := \dot{\bar\varphi}_{\kappa_m(t)}, 
\psi_t:=
\bar\varphi_{\kappa_m(t+a)-a}, \qquad
\dot\chi_t := 
\dot{\bar\varphi}_{\kappa_m(t+a)-a},
\qquad \chi_0=x.
$$
We have, with a unique constant $C = 2(KT+1)$ (see
(\ref{eqCa})) and for
any $\delta'$,
$$%\begin{equation}\label{eqC}
 P(\rho(X,\bar\varphi)<\delta') \le
 P(\rho(X^{\psi},\chi)<C\delta') \le
 P(\rho(X^{\psi,\Delta,a},\chi^{\Delta,a})<C\delta').
$$%\end{equation}
Denote $\delta''=C\delta'$. Let us choose $\delta''\le
\delta(\nu/20)$ (the notation from the Lemma  \ref{L-main}
is used), and
consider the following inequality, with the sequence $(\delta'_i,
\ 1\le i\le m)$, $\delta'_m = \delta''$, constructed via the value
$\nu/20$ instead of $\nu$ (compare to (\ref{modulus3}),
where the requirement related to $m_{\nabla H}$ could be
now dropped,
$$%\begin{equation}\label{eqC}
 P(\rho(X,\bar\varphi)<\delta'_1) \le
 E \prod_{i=1}^{m+1} 1(|X^{\psi,\Delta,a}_{t_k} -
\chi^{\Delta,a}_{t_k}| < \delta'_i), 
$$%\end{equation}
and $(t_k)$ are chosen as in the step 5. 
In particular, we require $4\delta'' = 4\delta'_{m+1} \le
\nu/20$, and $\sum_{i=1}^{m+1}\delta'_i\le
2\delta''$. Then,
due to the Lemma  \ref{L-main}
and using the same calculus as at the step 5, 
%%(\ref{it11}-\ref{it12}) (4-4??? not working!!!! nado
%%11-12), --
we get on the set
$\{|X^\psi_{t_{m}} - \chi_{t_{m}}| < \delta'_{m}\}$ and for
any $|\beta|\le b$,
\begin{eqnarray}\label{eq17a}
&E\left( I(|X^\psi_{t_{m+1}} -\chi_{t_{m+1}}| <
\delta'_{m+1}) \mid  {\cal F}_{t_{m}}\right) \nonumber
 \\ \nonumber \\
& \leq E^\beta \left( I(|X^\psi_{t_{m+1}} -\chi_{t_{m+1}}| <
\delta'_{m+1}) \exp\left( -\eps^{-2}\Delta_m \beta
\left((\chi_{t_{m+1}} - \chi_{t_{m}})/\Delta_m\right)
\nonumber \right. \right. \\ \nonumber \\
& \left. \left. - \eps^{-2}\Delta_m
(\tilde H(\psi_{t_{m}}, \psi_{t_{m}},\beta) - \nu/20) +
b\, \frac{\delta'_{m+1}+\delta'_{m}}{\varepsilon^{2}}\right) |
 {\cal F}_{t_{m}}\right)
\end{eqnarray}
(compare to (\ref{eq17})). The only change in comparison to
the step 5 is that now we want an upper bound, so  
indicators in the estimation will be just replaced by
$1$. Thus, we replace here 
$I(|X^\psi_{t_{m}} -\chi_{t_{m}}| < \delta'_{m+1})$ by $1$
and drop the expectation sign -- because there remains
nothing random in the expression -- then on the set
$\{|X^\psi_{t_{m}} - \chi_{t_{m}}| <
\delta'_{m}\}$ and for any $|\beta|\le b$ we get,
\begin{eqnarray}\label{eq17b}
&E\left(I(|X^\psi_{t_{m+1}} -\chi_{t_{m+1}}| <
\delta'_{m+1}) \mid 
 {\cal F}_{t_{m}}\right) \nonumber
 \\ \nonumber \\
& \leq \exp\left( -\eps^{-2}\Delta_m \beta
\left((\chi_{t_{m+1}}-\chi_{t_{m}})/\Delta_m \right)
\nonumber \right. \\ \nonumber \\
& \left. - \eps^{-2}\Delta_m (\tilde H(\psi_{t_{m}}, 
\psi_{t_{m}},\beta)) +
b\, \frac{\delta'_{m+1}+\delta'_{m}}{\varepsilon^{2}}\right),
\nonumber \\ \nonumber \\
& \leq \exp\left( -\eps^{-2}\Delta_m \beta
\left((\chi_{t_{m+1}} - \chi_{t_{m}})/\Delta_m \right)
\nonumber \right. \\ \nonumber \\
& \left. - \eps^{-2}\Delta_m (\tilde H(\psi_{t_{m}},
\psi_{t_{m}},\beta) - \nu/20) +
b\, \frac{\delta'_{m+1}+\delta'_{m}}{\varepsilon^{2}}\right),
\end{eqnarray}
once we have chosen $m_{\tilde H}(\delta'')\le \nu/20$ 
(remind that $m_{\tilde H} $
here means the modulus of continuity of the function 
$\tilde H(\cdot, \cdot, \beta)$ with respect to the first two 
variables on the set $|\beta|\le b+1$).

Let $\beta$ satisfy a condition,
\begin{eqnarray*} &\beta
(\chi_{t_{m+1}}-\chi_{t_{m}})/\Delta_m -
\tilde H(\psi_{t_{m}}, \psi_{t_{m}},\beta) 
 \\\\
&= \sup_{|\beta|\le b}\left(\beta
(\chi_{t_{m+1}}-\chi_{t_{m}})/\Delta_m -
\tilde H(\psi_{t_{m}},\psi_{t_{m}},\beta)\right)
 \\\\ &=
L^b(\psi_{t_{m}},\dot\chi_{t_{m} +}).
\end{eqnarray*}
Then, on the set $\{|X^\psi_{t_{m}} - \chi_{t_{m}}| <
\delta'_{m}\}$,
\begin{eqnarray}\label{eq17c}
&E\left(I(|X^\psi_{t_{m+1}} -\chi_{t_{m+1}}| <
\delta'_{m+1}) \mid
 {\cal F}_{t_{m}}\right) \nonumber
 \\ \nonumber \\
& \leq \exp\left( -\eps^{-2}\Delta_m
(L^b(\psi_{t_{m}},\dot\chi_{t_{m} +}) +\eps^{-2}\Delta_m
\frac{\nu}{20} +
b\, \frac{\delta'_{m+1} +
\delta'_{m}}{\varepsilon^{2}}\right).
\end{eqnarray}
Similarly and by induction and due to (\ref{elb}), we get
\begin{eqnarray}\label{eq17d}
&P(\rho(X,\chi)<\delta'_1) \nonumber
 \\ \nonumber \\
& \leq \exp \left(-\eps^{-2} \int\limits_0^T
%L^b(\bar\varphi_{\kappa_m(t)},\dot{\bar\varphi}_{
%\kappa_m(t) })\,dt
L^b(\psi_{t},\dot{\chi}_{t})\,dt
+ \eps^{-2}\nu/20 + 4b\,\varepsilon^{-2}\delta'_m \right)
\nonumber
 \\ \nonumber \\
&\le \exp\left( -\eps^{-2}(s  - \nu/5)\right).
\end{eqnarray}
This evidently contradicts (\ref{ge}).

~

%\item %20
\noindent
{\bf 17.}
Consider the case $\bar\varphi$ absolute continuous and
$S(\bar\varphi)=\infty$. In this case, due to monotone convergence
$L^b\to L$, there exist $b>0$, $m$ and $a\in [0,1]$ such
that
$$
\int_0^T L^b(\bar\varphi_t,\dot{\bar\varphi_t})\,dt \ge s -
\nu/20, \;  \int_0^T
L^b(\bar\varphi_{\kappa_m(t+a)-a},\dot{\bar\varphi}_{\kappa_m(t+a)-a})\,dt
\ge s - \nu/10.
$$
The rest is similar to the main case, $S(\bar\varphi)<\infty$, and
leads again to
$$
P(\rho(X,\bar\varphi)<\delta'_1)\le \exp\left( -\eps^{-2}(s  -
\nu/5)\right).
$$
This contradicts (\ref{ge}).

~

%\item %21
\noindent
{\bf 18.}
Consider the last possible case, $\bar\varphi$ not
absolute
continuous. In this case, for any constant $c$, in particular, for
$c=\|f\|_C+1$, there exist two values $0\le t_1 < t_2\le T$, such
that $|\bar\varphi_{t_2}-\bar\varphi_{t_1}|>c(t_2-t_1)$; indeed,
otherwise $\bar\varphi$ must be Lipschitz with
$|\dot{\bar\varphi}|\le c$. Therefore, for $\delta<(t_2-t_1)/2$,
probability $P(\rho(X,\bar\varphi)<\delta)$ necessarily
equals
zero, because the event $\{\rho(X,\bar\varphi)<\delta\}$ is empty.
This evidently contradicts (\ref{ge}). In all possible cases, we
got to contradictions. Hence, the assumption is wrong, that is,
the upper bound (\ref{eq4a}) holds true. The Theorem is proved.

%\end{enumerate}

~

%\vspace{1cm}

%\pagebreak
\begin{center}{\bf APPENDIX}\end{center}

~

{\bf A. Comments on the Lemma \ref{L-FW}}. To explain that the
Lemma \ref{L-FW} is valid without additional assumptions, we have
to review very briefly its proof and show those assumptions.

Let $0 = t_0 < t_0 < \ldots < t_m = T$ be a partition,
$\gamma_k(\beta) := \int_{t_{k-1}}^{t_{k}} H(\varphi_s,\beta)ds$,
$\ell_k(\alpha) = \sup_\beta (\alpha\beta - \gamma_k(\beta))$,
$A_k = \{\alpha: \; \ell_k(\alpha) < \infty\}$, $A_k^\circ$ its
interior with respect
to  the linear hull $L_{A_{k}}$.

The inequality $S(\varphi) = \int_0^T L(\varphi_t,\dot\varphi_t)dt
< \infty$ implies
$$
\sum\limits_{k=1}^m \sup\limits_\beta
\left((\varphi_{t_{k}}-\varphi_{t_{k-1}}) - \gamma_k(\beta)\right)
= \sum\limits_{k=1}^m \ell_k(\varphi_{t_{k}}-\varphi_{t_{k-1}})
\leq S(\varphi).
$$
Under additional assumption $A^\circ_k \not= \emptyset$ it is
proved in \cite{Freidlin and Wentzell (1984)} using the
arguments from
\cite{Rockafellar (1970)} that for any $\nu > 0$, there
exists a function
$\tilde \varphi$ such that $\rho(\varphi,\tilde\varphi) < \nu$ and
there exist $\beta_k$ such that

\begin{equation}\label{eq-FW1}
\ell_k(\tilde\varphi_{t_{k}} - \tilde\varphi_{t_{k-1}}) =
(\tilde\varphi_{t_{k}} - \tilde\varphi_{t_{k-1}})\beta_k -
\gamma_k(\beta_k)
\end{equation}
and
\begin{equation}\label{eq-FW21}
\tilde\varphi_{t_{k}}-\tilde\varphi_{t_{k-1}} =
\nabla\gamma_k(\beta_k).
\end{equation}
The proof goes well if $A^\circ_k \not= \emptyset \; \forall k$.

Let us show that the same is true if $A^\circ_k = \emptyset$ for
some $k$'s. The property $A^\circ_k = \emptyset$ is equivalent to
$\dim L_{A_{k}} = 0$. In this case, $\gamma_k(\beta) = c_k \beta$
with some $c_k\in R^d$. Hence, $\ell_k(\alpha_k) < \infty$ means
that $\ell_k(\alpha_k) = 0$ and for any other $\alpha$,
$\ell_k(\alpha) = + \infty$ and $\gamma_k(\beta) = \alpha_k\beta$.
So, we have
$$
\ell_k(\varphi_{t_{k}}-\varphi_{t_{k-1}}) = 0 =
(\varphi_{t_{k}}-\varphi_{t_{k-1}})\beta - \gamma_k(\beta)
$$
for {\it any} $\beta$. Let $\beta_k = 0$. Evidently,
$$
\varphi_{t_{k}}-\varphi_{t_{k-1}} = \nabla\gamma_k(\beta_k).
$$
Hence, in the case $A^\circ_k = \emptyset$, one just should
not change the curve $\varphi_s$ on the interval
$(t_{k-1},t_k)$; that is,
(\ref{eq-FW1}) %%??
and (\ref{eq-FW21}) are valid in this case also.

The rest of the proof remains unchanged. For any step
function
$\zeta$, one defines a piecewise linear $\chi$ by the formula
$$
\chi_0=\varphi_0, \quad \dot\chi_s = \nabla_\beta
H(\zeta_s,\beta_k)), \; t_{k-1}<s<t_k, \; k=1,2,\ldots ,m.
$$
Then it is shown that $\zeta^n \to \varphi$ implies $\chi^n\to
\varphi$ due to the property that the convergence of smooth convex
functions to the limit implies the convergence of their gradients.
Then there exists a partition such that this construction gives
one
$$
\int_0^T L(\zeta_t,\dot\chi_t)\,dt \leq S(\varphi)+\nu.
$$
So, the lemma holds true without additional assumptions. The
assertions about $\hat\zeta$ and $\hat\beta_s$ can be shown
similarly.

{\bf B. Comments on the property $A^\circ_k\not=\emptyset$, and
characterization of the set ${\cal L}^\circ[f,x]$}. Denote the
interior of $A(x) = \{\alpha: \, L(x,\alpha)<\infty\}$ with
respect to its linear hull $L_{A(x)}$ by $A^\circ(x)$. Then
$A^\circ_k=\emptyset \Longleftrightarrow
A^\circ(\varphi_{t_{k-1}})=\emptyset$. In this section we show the
following equivalence:
$$
card (f\in R^d:\; f=f(x,y), y\in M) = 1 \Longleftrightarrow \dim
L_{A(x)}=0 \Longleftrightarrow A^\circ(x)=\emptyset.
$$

Since $A(x)$ is convex, clearly the first two conditions are
equivalent.

If $\{f(x,\cdot)\}$ contains only one point then $H(x,\beta)$ is
linear with respect
to  $\beta$; hence, $A(x)$ consists of a unique point
and $A^\circ(x) = \emptyset$.

Now, let $\{f(x,\cdot)\}$ contain at least two different points,
say, $f(x,y_1) \not= f(x,y_2)$. Then there exists $1\leq k \leq d$
such that $(f(x,y_1)-f(x,y_2))_k \not= 0$. Denote $M_k = \sup_y
f^k(x,y), \; m_k = \inf_y f^k(x,y)$. Let $0 < \nu <
(f(x,y_1)-f(x,y_2))_k/2$. Take two points $y'$ and $y''$ such that
$f^k(x,y') < m_k + \nu/5$ and $f^k(x,y'') > M_k - \nu/5$. There
exist two open sets $B'\subset M$ and $B''\subset M$ such that
$\sup_{y\in B'} f^k(x,y) < m_k + \nu/4$ and $\inf_{y\in B''}
f^k(x,y) > M_k - \nu/4$.

Since the process $y^x_t$ is a nondegenerate ergodic diffusion,
there exists $\lambda >0$ such that
$$
P(y^x_s \in B',\; 1\leq s\leq t) \geq \lambda ^{t-1}, \quad
P(y^x_s \in B'',\; 1\leq s\leq t) \geq \lambda ^{t-1}, \quad t\to
\infty.
$$

Let $\beta = z \beta_k$ where $\beta_k\in E^d$ is a $k$th unit
coordinate vector and $z \in R$. Then for $z > 0$ we have,
\begin{eqnarray*}
z^{-1}t^{-1}\log E\exp(z\beta_k \int_0^t f(x,y_s^x)\,ds)
 \\ \\
\geq z^{-1}t^{-1}\log E\exp(z\beta_k \int_0^t f(x,y_s^x)\,ds)
I(y^x_s \in B'',\; 1\leq s\leq t)
 \\ \\
\geq z^{-1}t^{-1}\log \{\exp(z(M_k-\nu/2)t) \lambda ^{t-1}\} \\\\=
M_k-\nu/4 + \frac{t-1}{t} z^{-1}\log \lambda  \geq \frac{t-1}{t}
M_k - \nu/2,
\end{eqnarray*}
if $z$ is large enough. In other words, for large
positive $z$
one has $H(x,z\beta_k) \geq z (M_k-2\nu)$. Similarly, for large
negative $z$
\begin{eqnarray*}
|z|^{-1}t^{-1}\log E\exp(z\beta_k \int_0^t f(x,y_s^x)\,ds)
 \\ \\
\geq |z|^{-1}t^{-1}\log E\exp(z\beta_k \int_0^t f(x,y_s^x)\,ds)
I(y^x_s \in B'',\; 1\leq s\leq t)
 \\ \\
\geq |z|^{-1}t^{-1}\log \{\exp(z(m_k+\nu/4)t) \lambda ^{t-1}\}
 \\ \\
= -(m_k+\nu/4) + \frac{t-1}{t}|z|^{-1}\log \lambda  \geq -
\frac{t-1}{t}m_k - \nu/2,
\end{eqnarray*}
if $|z|$ is large enough. In other words, for negative $z$ with
large absolute values one has $H(x,z\beta_k) \geq z (m_k+\nu)$.
Therefore, $\{\alpha: \; \alpha = \beta_k \theta, \, m_k+\nu <
\theta < M_k - \nu\} \, \subset \, A(x)$.

On the other hand, it is obvious that if $\alpha = \beta_k
\theta$, $\theta\in R^1$, with $\theta > M_k$ or $\theta < m_k$,
then $L(x,\alpha)=\infty$, because $m_k z \le H(x,\beta_k z) \le
M_k z$, and, hence (say, if $\theta > M_k$), for $z >> 1$,
$$
\beta_k \theta \beta_k z - H(x,\beta_k z) \ge  (\theta - M_k) z
\to +\infty, \quad z\to +\infty.
$$

A similar calculus and similar inequalities are valid for
any unit vector
$\beta_0$. This shows, in particular, that $\dim L_A(x) = \dim
L_f(x)$, and, moreover, that $L_A(x) = L_f(x)$. Since $A(x)$ is
convex, it shows also that the interior $A^\circ(x)$ with
respect
to 
$L_{A(x)}$ is not empty, except for only the case
dim$(L_{A(x)}) =
1$. Hence, the third condition is equivalent to the second one and
to the first.
 
So, the condition $A^\circ_k \not= \emptyset$ is always satisfied
if the set $\{f(x,\cdot)\}$ for any $x$ consists of more than one
point. In fact, if $card \{f(x,\cdot)\}=1$ for {\it any} $\,x$
then $f$ does not depend on $y$. In this case, one has nothing to
average.

Notice that our considerations above provide the following
description of the set ${\cal L}^\circ[f,x]$:
\begin{eqnarray*}
&{\cal L}^\circ[f,x] = \{\alpha\in R^d: \, m_\beta(x) < \langle
\alpha,\beta\rangle <M_\beta(x), \; \forall %\beta\in R^d,
|\beta|=1, \; \mbox{with}\; m_\beta(x) < M_\beta(x), \\
&\mbox{and}\; \langle \alpha,\beta\rangle = M_\beta(x), \; \forall
|\beta|=1, \; \mbox{with}\; m_\beta(x) = M_\beta(x) \},
\end{eqnarray*}

\vspace{2mm}

\noindent where $m_\beta(x) := \inf\limits_y \langle
\frac{\beta}{|\beta|}, f(x,y)\rangle, \, M_\beta(x) :=
\sup\limits_y \langle \frac{\beta}{|\beta|}, f(x,y)\rangle$.
Moreover, it can be shown similarly that for any $x, \tilde x$
(although we do not need it here),
$$
{\cal L}^\circ[f,x,\tilde x] = {\cal L}^\circ[f,x].
$$

~

%%%%%%%%%%%%%%%%%%%%%%%%

%%%%%%%%%%%%%%%%%%%%%%%%

{\bf C. About $\hat\alpha_s \in {\cal L}^\circ[f,\varphi_s]$}. Let
$x=\varphi_s$, $\hat\alpha = \hat\alpha[x,\dot\chi]$ as described
in the proof of the theorem 1. If we show that for any direction
$v$ (a unit vector) satisfying the property $m_v < M_v$,
%$$\inf (vf(x,y),\, y\in R^d) =: m_v <
%M_v:= \sup (vf(x,y),\, y\in R^d),
%$$ %for any $\nu>0$
the strict double inequality holds true
$$
m_v  < \partial H(x,zv)/\partial z|_{z=0} < M_v,
$$
$z\in R^1$, then it would follow $\hat\alpha_s \in {\cal
L}^\circ[f,\varphi_s]$. Let $\nu>0$ and again two open sets $B'$
and $B''$ be chosen such that $\sup_{y\in B'} vf(x,y) < m_v +
\nu/2$, and $\inf_{y\in B''} vf(x,y) > M_v- \nu/2$. Let
$\mu_{inv}(B'')$ be invariant measure for the event $\{y^x_t\in
B''\}$. We can choose $\nu$ and correspondingly $B''$ so that
$\mu_{inv}(B'') < 1$. Then, due to large deviation asymptotics for
the process $y^x_t$, for any $\mu_{inv}(B'') < \zeta < 1$ there
exists $\lambda>0$ such that
$$
P\left(t^{-1}\int_0^t 1(y^x_s\in  B'')\,ds \ge \zeta\right) \le
\exp(-\lambda t), \quad t\ge t_{\zeta}.
$$
Denote $A_\zeta = \left\{t^{-1}\int_0^t 1(y^x_s\in  B'')\,ds
<\zeta\right\}$, $A^c_\zeta = \left\{t^{-1}\int_0^t 1(y^x_s\in
B'')\,ds \ge \zeta\right\}$, then for $z>0$,
\begin{eqnarray*}
E\exp(z v \int_0^t f(x,y^x)\,ds)
 \\ \\
\le E\exp(z \int_0^t
 \left(M_v 1(y^x_s\in B'') + (M_v-\nu) 1(y^x_s\not\in B'')\right)
 \,ds)\,1(A^c_\zeta)
 \\ \\
+E\exp(z \int_0^t
 \left(M_v 1(y^x_s\in B'') + (M_v-\nu) 1(y^x_s\not\in B'')\right)
 \,ds)\, 1(A_\zeta)
  \\ \\
\le E\exp(z t M_v + z t (M_v-\nu)) \,1(A^c_\zeta)
 %\\ \\
+E\exp(z t M_v \zeta  + zt(M_v-\nu))\, 1(A_\zeta)
 \\ \\
\le \exp(z t M_v + z t (M_v-\nu) -z\lambda t/z)
 %\\ \\
+\exp(z t M_v \zeta  + zt(M_v-\nu)),
\end{eqnarray*}
hence,
$$
\limsup_{z\to 0}\,\limsup_{t\to\infty}\, (tz)^{-1}
\ln E\exp(z v
\int_0^t f(x,y^x)\,ds) < M_v.
$$
Similarly, using $B'$ one can get
$$
\liminf_{z\to 0}\,\liminf_{t\to\infty} \,(tz)^{-1} \ln
 E\exp(z v
\int_0^t f(x,y^x)\,ds) > m_v.
$$
Thus,
$$
m_v < \partial H(x,zv)/\partial z|_{z=0} < M_v.
$$
Therefore, $\hat\alpha \in {\cal L}^\circ[x,f]$.

~

%\pagebreak

\section*{Acknowledgements}
The author is grateful to Professor Yuri Kifer who
initiated the process of correction and to the
unknown referee of this verson of the paper for very helpful
remarks that resulted in the final simplifications and
improvements.

%\begin{center}{REREFERCES}\end{center}

\end{document}